\documentclass[11pt]{amsart}
\usepackage{xypic}

\newtheorem{proposition}{Proposition}[section]
\newtheorem{lemma}[proposition]{Lemma}

\newtheorem{theorem}[proposition]{Theorem}

\theoremstyle{definition}
\newtheorem{definition}[proposition]{Definition}
\newtheorem{example}[proposition]{Example}

\theoremstyle{remark}
\newtheorem{remark}[proposition]{Remark}

\newcommand{\thlabel}[1]{\label{th:#1}}
\newcommand{\thref}[1]{Theorem~\ref{th:#1}}
\newcommand{\selabel}[1]{\label{se:#1}}
\newcommand{\seref}[1]{Section~\ref{se:#1}}
\newcommand{\lelabel}[1]{\label{le:#1}}
\newcommand{\leref}[1]{Lemma~\ref{le:#1}}
\newcommand{\prlabel}[1]{\label{pr:#1}}
\newcommand{\prref}[1]{Proposition~\ref{pr:#1}}

\newcommand{\exlabel}[1]{\label{ex:#1}}
\newcommand{\exref}[1]{Example~\ref{ex:#1}}
\newcommand{\delabel}[1]{\label{de:#1}}

\newcommand{\eqlabel}[1]{\label{eq:#1}}
\newcommand{\equref}[1]{(\ref{eq:#1})}

\def\mapright#1{\smash{\mathop{\longrightarrow}\limits^{#1}}}

\def\equal#1{\smash{\mathop{=}\limits^{#1}}}

\newcommand{\can}{{\rm can}}
\newcommand{\Hom}{{\rm Hom}}
\newcommand{\HOM}{{\rm HOM}}

\newcommand{\END}{{\rm END}}

\def\lan{\langle}
\def\ran{\rangle}
\def\ot{\otimes}

\def\GG{{\mathbb G}}

\def\MM{{\mathbb M}}
\def\rightact{\hbox{$\leftharpoonup$}}

\newcommand{\Cc}{\mathcal{C}}
\newcommand{\Dd}{\mathcal{D}}

\newcommand{\Mm}{\mathcal{M}}

\def\*C{{}^*\hspace*{-1pt}{\Cc}}

\def\text#1{{\rm {\rm #1}}}

\def\ul{\underline}

\begin{document}
\title[Group corings]{Group corings}
\author{S. Caenepeel}
\address{Faculty of Engineering,
Vrije Universiteit Brussel, B-1050 Brussels, Belgium}
\email{scaenepe@vub.ac.be}
\urladdr{http://homepages.vub.ac.be/\~{}scaenepe/}
\author{K. Janssen}
\address{Faculty of Engineering,
Vrije Universiteit Brussel, B-1050 Brussels, Belgium}
\email{krjansse@vub.ac.be}
\urladdr{http://homepages.vub.ac.be/\~{}krjansse/}
\author{S.H. Wang}
\address{Department of mathematics, Southeast University,
 Nanjing 210096, China}
\email{shuanhwang2002@yahoo.com}



\thanks{This research was supported by the research project G.0622.06 ``Deformation quantization methods
for algebras and categories with applications to quantum mechanics" from
FWO-Vlaanderen. The third author was partially supported by
the SRF (20060286006) and the FNS (10571026)}

\begin{abstract}
We introduce group corings, and study functors between categories of comodules over
group corings, and the relationship to graded modules over graded rings. Galois group
corings are defined, and a Structure Theorem for the $G$-comodules over a Galois
group coring is given. We study (graded) Morita contexts associated to a group coring.
Our theory is applied to group corings associated to a comodule algebra over a
Hopf group coalgebra.
\end{abstract}

\maketitle

\section*{Introduction}\selabel{0}
Group coalgebras and Hopf group coalgebras were introduced by Turaev \cite{Turaev}.
A systematic algebraic study of these new structures has been carried out in recent
papers by Virelizier, Zunino, and the third author (see
for example \cite{Virelizier,Wang1,Wang2,Wang3,Zunino1,Zunino2}). Many results from classical
Hopf algebra theory can be generalized to Hopf group coalgebras; this has been
explained in a paper by the first author and De Lombaerde \cite{CDL}, where it was
shown that Hopf group coalgebras are in fact Hopf algebras in a suitable symmetric
monoidal category.\\
In \cite{Wang4}, the third author investigated how Hopf-Galois theory can be developed in
the framework of Hopf group coalgebras. A definition of Hopf-Galois extension was
presented; the requirement is that a set of canonical maps, indexed by the elements of
the underlying group, has to be bijective. 
One aspect in the present theory that is not satisfactory
is the lack of an appropriate Structure Theorem: an important result in Hopf-Galois
theory states that the category of relative Hopf modules over a faithfully flat Hopf-Galois
extension is equivalent to the category of modules over the coinvariants. So far,
no such result is known in the framework of Hopf group coalgebras.\\
Corings were introduced by Sweedler \cite{Sweedler65}, and were revived
recently by Brzezi\'nski \cite{Br3}. One of the important observations is that
coring theory provides an elegant 
approach to descent theory and Hopf-Galois theory (see \cite{Br3,BrW,C03}).
The aim of this paper is to develop Galois theory for group corings, and to apply
it to Hopf group coalgebras.\\
A $G$-$A$-coring (or group coring) consists of a set of $A$-bimodules indexed by a group $G$,
together with a counit map, and a set of diagonal maps indexed by $G\times G$,
with appropriate axioms (see \seref{1}).
A first remarkable observation is the fact that we can introduce
two dif\emph{}ferent types of comodules over a group coring $\ul{\Cc}$. $\ul{\Cc}$-comodules
consist of a single $A$-module, with a set of structure maps indexed by $G$, while
$G$-$\ul{\Cc}$-comodules consist of a set of $A$-modules indexed by $G$, together with structure maps
indexed by $G\times G$. We have a pair of adjoint functors between the two categories of
comodules (see \prref{1.2}). This remarkable fact can be explained by duality arguments.
Dualizing the definition of a $G$-$A$-coring, we obtain $G$-$A$-rings; in the
coalgebra case, this was observed in \cite{Zunino1}. In contrast to $G$-corings,
$G$-rings are a well-known concept: in fact there is a categorical correspondence
between $G$-$A$-rings and $G$-graded $A$-rings. The graded ring corresponding
to a $G$-ring is the so-called ``packed form" of the $G$-ring, in the terminology of
\cite{Zunino1}. We don't have a similar correspondence between group corings and
graded corings, unless the group $G$ in question is finite. This indicates that duality
properties for module categories over graded rings have to be studied from the point of
view of group corings, rather than graded corings.\\
Over a graded ring, one can study ordinary modules and graded modules, and there
exists an adjoint pair between the two categories. We also have functors from the
categories of $\ul{\Cc}$-comodules (resp.  $G$-$\ul{\Cc}$-comodules) to
modules (resp. graded modules) over the dual graded ring of $\ul{\Cc}$. All these
functors appear in a commutative diagram of functors (see \prref{4.6}). The functor
between the category of $G$-$\ul{\Cc}$-comodules and graded modules is an equivalence
if every $\Cc_\alpha$ (or equivalently, every homogeneous part of the dual graded
ring) is finitely generated and projective as an $A$-module (see \prref{4.5}).
These properties of (co)module categories are studied in Sections \ref{se:1}, \ref{se:3}
and \ref{se:4}.\\
An important class of group corings, called cofree group corings, is investigated in \seref{2}.
Basically, these are corings for which all the underlying bimodules are isomorphic.
A cofree coring is - up to isomorphism - determined by $\Cc_e$, its part of degree $e$;
its left dual is the group ring ${}^*\Cc_e[G]$ over the left dual of $\Cc_e$
(\prref{4.7}). The category of $G$-$\Cc$-comodules is equivalent to the category of
comodules over $\Cc_e$ (\thref{2.2}). This is an analog of the well-known fact that,
for a group ring $R[G]$, the category of graded $R[G]$-modules is equivalent to the
category of $R$-modules.\\
In \seref{5}, we introduce the notion of grouplike family of a group coring. Grouplike families
correspond bijectively to $\ul{\Cc}$-comodule structures on $A$. Fixing a grouplike
family, we can introduce the coinvariant subring $T$ of $A$. We have two pairs
of adjoint functors, one connecting modules over the coinvariants to right 
$\ul{\Cc}$-comodules (\prref{5.3}), and another one connecting modules over the coinvariants to right 
$G$-$\ul{\Cc}$-comodules (\prref{5.4}). It can be established when the latter adjoint pair,
denoted $(F_7,G_7)$,
is a pair of inverse equivalences. Given a group coring $\ul{\Cc}$ with a fixed grouplike
family, we can define a canonical morphism of group corings between the cofree coring
built on the Sweedler canonical coring and the coring $\ul{\Cc}$. If $F_7$ is an
equivalence, then this canonical morphism is an isomorphism (\prref{5.7}). In this case, we call
our group coring a Galois group coring. This is equivalent to $\ul{\Cc}$ being cofree,
and $\Cc_e$ being a Galois coring. The Structure Theorem \ref{th:5.12} is our
main result.
Basically, if $\ul{\Cc}$ is Galois, and $A$ is faithfully flat over the coinvariants, then $F_7$
is an equivalence.\\
Morita theory plays an important role in this theory. To a group coring with a fixed grouplike
family, we can associate several Morita contexts. Two of them are classical Morita contexts,
and have been studied in a special situation (see \seref{10}) in \cite{Wang4}. But the
natural Morita contexts are in fact graded Morita contexts. In \seref{6}, we give some
generalities on graded Morita contexts; in Sections \ref{se:7} and \ref{se:8}, we discuss the Morita contexts and their relationship.\\
In some situations, the Galois property of a group coring can be characterized by the
graded Morita contexts associated to it. We study this in \seref{9}. In the final \seref{10},
we briefly discuss the situation where $\ul{\Cc}$ is a group coring $A\ot \ul{H}$ associated to a
right $\ul{H}$-comodule algebra $A$, where $\ul{H}$ is a Hopf group coalgebra, as
introduced in \cite{Turaev}. We show that this group coring is Galois if and only if
$A$ is an $\ul{H}$-Galois extension in the sense of \cite{Wang4}. This entails a Structure
Theorem for relative group Hopf modules; we also describe the dual of the group coring
$A\ot \ul{H}$.\\

Throughout this paper, we will adopt the following notational conventions.\\
For an object $M$ in a category, $M$ will also denote the identity morphism on $M$.\\
Let $G$ be a group and $M$ a (right) $A$-module. We will often need collections of
$A$-modules isomorphic to $M$ and indexed by $G$. We will consider these modules
as isomorphic, but distinct. Let $M\times\{\alpha\}$ be the module with index $\alpha$.
We then have isomorphisms
$$\mu_\alpha:\ M\to M\times\{\alpha\},~~\mu_\alpha(m)=(m,\alpha).$$
We can then write $M\times\{\alpha\}=\mu_\alpha(M)$. $\mu$ can be considered as a
dummy variable, we will also use the symbols $\gamma,\nu,\ldots$. We will identify
$M$ and $M\times \{e\}$ using $\mu_e$.

\section{Group corings and comodules}\selabel{1}
Let $G$ be a group, and $A$ a ring with unit. The unit element of $G$ will be denoted
by $e$. A $G$-group $A$-coring (or shortly a $G$-$A$-coring) $\ul{\Cc}$ is a family $(\Cc_\alpha)_{\alpha\in G}$ of $A$-bimodules
together with a family of bimodule maps
$$\Delta_{\alpha,\beta}:\ \Cc_{\alpha\beta}\to \Cc_\alpha\ot_A \Cc_\beta~~;~~
\varepsilon:\ \Cc_e\to A,$$
such that
\begin{equation}\eqlabel{1.1.1}
(\Delta_{\alpha,\beta}\ot_A\Cc_\gamma)\circ \Delta_{\alpha\beta,\gamma}=
(\Cc_\alpha\ot_A \Delta_{\beta,\gamma})\circ \Delta_{\alpha,\beta\gamma}
\end{equation}
and
\begin{equation}\eqlabel{1.1.2}
(\Cc_\alpha\ot_A\varepsilon)\circ\Delta_{\alpha,e}=\Cc_\alpha=(\varepsilon\ot_A \Cc_\alpha)
\circ \Delta_{e,\alpha},
\end{equation}
for all $\alpha,\beta,\gamma\in G$.
We use the following Sweedler-type notation for the comultiplication maps $\Delta_{\alpha,\beta}$:
$$\Delta_{\alpha,\beta}(c)=c_{(1,\alpha)}\ot_A c_{(2,\beta)},$$
for all $c\in \Cc_{\alpha\beta}$. Then \equref{1.1.2} takes the form
\begin{equation}\eqlabel{1.1.3}
c_{(1,\alpha)}\varepsilon(c_{(2,e)})=c=\varepsilon(c_{(1,e)})c_{(2,\alpha)}.
\end{equation}
\equref{1.1.1} justifies the following notation:
$$
((\Delta_{\alpha,\beta}\ot_A\Cc_\gamma)\circ \Delta_{\alpha\beta,\gamma})(c)=
((\Cc_\alpha\ot_A \Delta_{\beta,\gamma})\circ \Delta_{\alpha,\beta\gamma})(c)=
c_{(1,\alpha)}\ot_A c_{(2,\beta)}\ot_A c_{(3,\gamma)},$$
for all $c\in \Cc_{\alpha\beta\gamma}$. If $\ul{\Cc}$ is a $G$-$A$-coring,
then $\Cc=\Cc_e$ is an $A$-coring, with comultiplication $\Delta_{e,e}$ and counit
$\varepsilon$.\\
A morphism between two $G$-$A$-corings $\ul{\Cc}$ and $\ul{\Dd}$
consists of a family of $A$-bimodule maps $(f_\alpha)_{\alpha\in G}$,
$f_\alpha:\ \Cc_\alpha\to \Dd_\alpha$ such that
$$(f_\alpha\ot_A f_\beta)\circ\Delta_{\alpha,\beta}=\Delta_{\alpha,\beta}\circ f_{\alpha\beta}
~~{\rm and}~~\varepsilon\circ f_e=\varepsilon.$$

If $A$ is a commutative ring, and $ac=ca$, for all $\alpha\in G$, $a\in A$ and $c\in \Cc_\alpha$,
then $\ul{\Cc}$ is called a $G$-coalgebra, cf. \cite{Turaev}.\\

Over a group coring, we can define two dif\emph{}ferent types of comodules. A right $\ul{\Cc}$-comodule
is a right $A$-module $M$ together with a family of right $A$-linear maps
$(\rho_\alpha)_{\alpha\in G}$, $\rho_\alpha:\ M\to M\ot_A \Cc_\alpha$, such that
\begin{equation}\eqlabel{1.1.4}
(M\ot_A \Delta_{\alpha,\beta})\circ\rho_{\alpha\beta}=
(\rho_\alpha\ot_A\Cc_\beta)\circ\rho_\beta
\end{equation}
and
\begin{equation}\eqlabel{1.1.5}
(M\ot_A\varepsilon)\circ\rho_e=M.
\end{equation}
We use the following Sweedler-type notation:
$$\rho_\alpha(m)=m_{[0]}\ot_A m_{[1,\alpha]}.$$
\equref{1.1.4} justifies the notation
$$((M\ot_A \Delta_{\alpha,\beta})\circ\rho_{\alpha\beta})(m)=
((\rho_\alpha\ot_A\Cc_\beta)\circ\rho_\beta)(m)=m_{[0]}\ot_A m_{[1,\alpha]}\ot_A m_{[2,\beta]},$$
and \equref{1.1.5} is equivalent to $m_{[0]}\varepsilon(m_{[1,e]})=m$, for all $m\in M$.\\
A morphism of right $\ul{\Cc}$-comodules is a right $A$-linear map $f:\ M\to N$ satisfying the
condition
\begin{equation}\eqlabel{1.1.6}
(f\ot_A\Cc_\alpha)\circ\rho_\alpha=\rho_\alpha\circ f,
\end{equation}
for all $\alpha\in G$. $\Mm^{\ul{\Cc}}$ will be our notation for the category of
right $\ul{\Cc}$-comodules.\\

A right $G$-$\ul{\Cc}$-comodule $\ul{M}$ is a family of right $A$-modules
$(M_\alpha)_{\alpha\in G}$, together with a family of right $A$-linear maps
$$\rho_{\alpha,\beta}:\ M_{\alpha\beta}\to M_\alpha\ot_A \Cc_\beta$$
such that
\begin{equation}\eqlabel{1.1.7}
(M_\alpha\ot_A\Delta_{\beta,\gamma})\circ \rho_{\alpha,\beta\gamma}=
(\rho_{\alpha,\beta}\ot_A\Cc_\gamma)\circ\rho_{\alpha\beta,\gamma}
\end{equation}
and
\begin{equation}\eqlabel{1.1.8}
(M_\alpha\ot_A\varepsilon)\circ \rho_{\alpha,e}=M_\alpha
\end{equation}
for all $\alpha,\beta,\gamma\in G$. We now use the following Sweedler-type notation:
$$\rho_{\alpha,\beta}(m)=m_{[0,\alpha]}\ot_A m_{[1,\beta]},$$
for $m\in M_{\alpha\beta}$. \equref{1.1.7} justifies the notation
\begin{eqnarray*}
&&\hspace*{-2cm}
((M_\alpha\ot_A\Delta_{\beta,\gamma})\circ \rho_{\alpha,\beta\gamma})(m)=
((\rho_{\alpha,\beta}\ot_A\Cc_\gamma)\circ\rho_{\alpha\beta,\gamma})(m)\\
&=&
m_{[0,\alpha]}\ot_A m_{[1,\beta]}\ot_A m_{[2,\gamma]},
\end{eqnarray*}
for $m\in M_{\alpha\beta\gamma}$. \equref{1.1.8} implies that
$m_{[0,\alpha]}\varepsilon(m_{[1,e]})=m$, for all $m\in M_\alpha$. A morphism between two
right $G$-$\ul{\Cc}$-comodules $\ul{M}$ and $\ul{N}$ is a family of right $A$-linear
maps $f_\alpha:\ M_\alpha\to N_\alpha$ such that
$$(f_\alpha\ot_A\Cc_\beta)\circ \rho_{\alpha,\beta}=\rho_{\alpha,\beta}\circ
f_{\alpha\beta}.$$
The category of right $G$-$\ul{\Cc}$-comodules will be denoted by $\Mm^{G,\ul{\Cc}}$. 

\begin{proposition}\prlabel{1.2}
We have a pair of adjoint functors $(F_1,G_1)$ between the categories
$\Mm^{G,\ul{\Cc}}$ and $ \Mm^{\ul{\Cc}}$. Moreover, if $G$ is a finite group, then $(F_1,G_1)$ is a Frobenius pair of functors, i.e. $F_1$ is also a right adjoint of $G_1$.
\end{proposition}

\begin{proof}
Take $\ul{M}=(M_\alpha)_{\alpha\in G}\in \Mm^{G,\ul{\Cc}}$, and define
$$F_1(\ul{M})=\bigoplus_{\alpha\in G} M_\alpha=M.$$
The coaction maps $\rho_\alpha:\ M\to M\ot_A\Cc_\alpha$ are defined as follows:
for $m\in M_\beta$, let
\begin{equation}\eqlabel{1.2.1}
\rho_\alpha(m)=m_{[0,\beta\alpha^{-1}]}\ot_A m_{[1,\alpha]}.
\end{equation}
Otherwise stated,
$\rho_\alpha=\bigoplus_{\beta\in G} \rho_{\beta\alpha^{-1},\alpha}$.
Let us show that (\ref{eq:1.1.4},\ref{eq:1.1.5}) hold. For all $m\in M_\gamma$,
we compute that
\begin{eqnarray*}
&&\hspace*{-2cm}
((\rho_\alpha\ot_A\Cc_\beta)\circ\rho_\beta)(m)=
(\rho_\alpha\ot_A\Cc_\beta)(m_{[0,\gamma\beta^{-1}]}\ot_A m_{[1,\beta]})\\
&=& m_{[0,\gamma\beta^{-1}\alpha^{-1}]}\ot_A m_{[1,\alpha]}\ot_A m_{[2,\beta]}\\
&=&((M\ot_A\Delta_{\alpha,\beta})\circ\rho_{\alpha\beta})(m);\\
&&\hspace*{-2cm}
((M\ot_A\varepsilon)\circ\rho_e)(m)=m_{[0,\gamma]}\varepsilon(m_{[1,e]})=m.
\end{eqnarray*}
For a morphism $\ul{f}:\ \ul{M}\to \ul{N}$ in $\Mm^{G,\ul{\Cc}}$, we simply define
$$F_1(\ul{f})=\bigoplus_{\alpha\in G} f_\alpha.$$
Let us now define $G_1$. For $M\in \Mm^{\ul{\Cc}}$, let
$G_1(M)_\alpha=\mu_\alpha(M)$, where we use the notation introduced at the end of
the introduction. The coaction maps $\rho_{\alpha,\beta}:\ \mu_{\alpha\beta}(M)\to
 \mu_{\alpha}(M)\ot_A \Cc_\beta$ are defined by
 \begin{equation}\eqlabel{1.2.2}
 \rho_{\alpha,\beta}(\mu_{\alpha\beta}(m))=\mu_\alpha(m_{[0]})\ot_A m_{[1,\beta]},
 \end{equation}
for all $m\in M$. The formulas (\ref{eq:1.1.7},\ref{eq:1.1.8}) hold since
\begin{eqnarray*}
&&\hspace*{-2cm}
((M_\alpha\ot_A \Delta_{\beta,\gamma})\circ \rho_{\alpha,\beta\gamma})(\mu_{\alpha\beta\gamma}
(m))\\
&=& (M_\alpha\ot_A \Delta_{\beta,\gamma})(\mu_\alpha(m_{[0]})\ot_A m_{[1,\beta\gamma]})\\
&=& \mu_\alpha(m_{[0]})\ot_A m_{[1,\beta]}\ot_A m_{[2,\gamma]}\\
&=& \rho_{\alpha,\beta}(\mu_{\alpha\beta}(m_{[0]}))\ot_A m_{[1,\gamma]}\\
&=& ((\rho_{\alpha,\beta}\ot_A \Cc_\gamma)\circ \rho_{\alpha\beta,\gamma})
(\mu_{\alpha\beta\gamma}(m));\\
&&\hspace*{-2cm}
((M_\alpha\ot_A\varepsilon)\circ\rho_{\alpha,e})(\mu_\alpha(m))\\
&=& (M_\alpha\ot_A\varepsilon)(\mu_\alpha(m_{[0]})\ot_A m_{[1,e]})\\
&=& \mu_\alpha(m_{[0]})\varepsilon(m_{[1,e]})=\mu_\alpha(m_{[0]}\varepsilon(m_{[1,e]}))=
\mu_\alpha(m).
\end{eqnarray*}
On the morphisms, $G_1$ is defined as follows: for $f:\ M\to N$ in $\Mm^{\ul{\Cc}}$, we put 
$$G_1(f)=(\nu_\alpha\circ f \circ \mu_\alpha^{-1})_{\alpha\in G}.$$
Take $\ul{M}\in \Mm^{G,\ul{\Cc}}$ and $N\in \Mm^{\ul{\Cc}}$, and consider the map
$$\psi:\ \Hom^{\ul{\Cc}}(F_1(\ul{M}),N)\to \Hom^{G,\ul{\Cc}}(\ul{M},G_1(N))$$
defined as follows. For $f:\ \bigoplus_{\alpha\in G} M_\alpha\to N$, let
$$\psi(f)_\alpha=\nu_\alpha\circ f\circ i_\alpha :\ M_\alpha\to G_1(N)_\alpha=\nu_\alpha(N),$$
where $i_\alpha: M_\alpha\to \bigoplus_{\alpha\in G} M_\alpha$ is the canonical injection.
Now consider the map
$$\phi:\ \Hom^{G,\ul{\Cc}}(\ul{M},G_1(N))\to \Hom^{\ul{\Cc}}(F_1(\ul{M}),N),$$
defined as follows: for $\ul{g}=(g_\alpha)_{\alpha\in G}:\ \ul{M}\to G_1(N)$,
let
$$\phi(\ul{g})(m)=\sum_{\alpha\in G} (\nu_{\alpha}^{-1}\circ g_{\alpha}\circ p_\alpha)(m),$$
where now $p_\alpha:\ \bigoplus_{\alpha\in G} M_\alpha\to M_\alpha$ is the canonical
projection. Straightforward computations show that $\psi$ and $\phi$ are
well-defined. They are inverses, since
$$
\phi(\psi(f))(m)=\sum_{\alpha\in G} (\nu_\alpha^{-1}\circ \nu_\alpha\circ f\circ i_\alpha\circ p_\alpha)(m)
= f(\sum_{\alpha\in G} (i_\alpha\circ p_\alpha)(m))=f(m),
$$
for all $m\in \bigoplus_{\alpha\in G} M_\alpha$, and
\begin{eqnarray*}
&&\hspace*{-15mm}
\psi(\phi(\ul{g}))_\alpha(m)=(\nu_\alpha\circ \phi(\ul{g})\circ i_\alpha)(m)\\
&=& \sum_{\beta\in G} (\nu_\alpha\circ \nu_\beta^{-1}\circ g_\beta\circ p_\beta\circ i_\alpha)(m)
= (\nu_\alpha\circ \nu_\alpha^{-1}\circ g_\alpha)(m)=g_\alpha(m),
\end{eqnarray*}
for all $\alpha\in G$ and $m\in M_\alpha$. It is easy to show that $\psi$ and
$\phi$ define natural transformations.
Let us describe the unit $\eta_1$ and the counit $\varepsilon_1$ of the adjunction.
For $\ul{M}\in \Mm^{G,\ul{\Cc}}$, we have
$$\eta_{1,\ul{M},\beta}=\mu_\beta\circ i_\beta:\ 
M_\beta\to \mu_{\beta}(\bigoplus_{\alpha\in G} M_\alpha);$$
for $N\in \Mm^{\ul{\Cc}}$, we have
$$\varepsilon_{1,N}=\sum_{\alpha\in G} \mu_{\alpha}^{-1}\circ p_\alpha:\ \bigoplus_{\alpha\in G}\mu_\alpha(N)\to N.$$
To prove the final statement, let us assume that $G$ is finite. Take $\ul{M}\in \Mm^{G,\ul{\Cc}}$ and $N\in \Mm^{\ul{\Cc}}$, and consider the map
$$\Phi:\ \Hom^{G,\ul{\Cc}}(G_1(N),\ul{M})\to \Hom^{\ul{\Cc}}(N,F_1(\ul{M})),$$
defined by
$$\Phi(\ul{g})(n)=\sum_{\alpha\in G}(i_\alpha \circ g_\alpha\circ \mu_\alpha)(n)\in \bigoplus_{\alpha\in G}M_\alpha,$$
for all morphisms $\ul{g}=(g_\alpha)_{\alpha\in G}:\ G_1(N)\to \ul{M}$ in $\Mm^{G,\ul{\Cc}}$.
Now consider the map
$$\Psi:\ \Hom^{\ul{\Cc}}(N,F_1(\ul{M}))\to \Hom^{G,\ul{\Cc}}(G_1(N),\ul{M}),$$
defined by 
$$\Psi(f)_\alpha=p_\alpha\circ f \circ \mu_\alpha^{-1}:\mu_\alpha(N)\to M_\alpha,$$
for all right $\ul{\Cc}$-colinear maps $f:N\to \bigoplus_{\alpha\in G}M_\alpha.$ One can check that $\Phi$ and $\Psi$ are well-defined. Let us check that $\Phi$ and $\Psi$ are mutually inverse:
\begin{eqnarray*}
&&\hspace*{-5mm}
\Phi(\Psi(f))(n)=\sum_{\alpha\in G} (i_\alpha\circ \Psi(f)_\alpha \circ \mu_\alpha)(n)\\
&=& \sum_{\alpha\in G} (i_\alpha\circ p_\alpha \circ f \circ \mu_\alpha^{-1} \circ \mu_\alpha)(n)
= \sum_{\alpha\in G} (i_\alpha\circ p_\alpha)(f(n))=f(n),
\end{eqnarray*}
for all $n\in N$, and, for all $\alpha \in G$,
\begin{eqnarray*}
&&\hspace*{-2cm}
\Psi(\Phi(\ul{g}))_\alpha=p_\alpha \circ \Phi(\ul{g})\circ \mu_\alpha^{-1}= \sum_{\beta \in G} (p_\alpha \circ i_\beta \circ g_\beta \circ \mu_\beta \circ \mu_\alpha^{-1})\\
&=& g_\alpha \circ \mu_\alpha \circ \mu_\alpha^{-1}= g_\alpha.
\end{eqnarray*}
Let us finally describe the unit $\nu_1$ and the counit $\zeta_1$ of this adjunction.
For $N\in \Mm^{\ul{\Cc}}$, we have
$$\nu_{1,N}=\sum_{\alpha\in G} i_{\alpha}\circ \mu_\alpha:\ N\to \bigoplus_{\alpha \in G}\mu_\alpha(N);$$
for $\ul{M}\in \Mm^{G,\ul{\Cc}}$, we have
$$\zeta_{1,\ul{M},\beta}=p_\beta\circ \mu_\beta^{-1}:\ \mu_\beta(\bigoplus_{\alpha \in G}M_\alpha)\to M_\beta. $$
\end{proof}

\section{Cofree group corings}\selabel{2}
\begin{definition}\delabel{2.1}
A $G$-$A$-coring $\ul{\Cc}$ is called cofree if there exist $A$-bimodule isomorphisms
$\gamma_\alpha:\ \Cc=\Cc_e\to \Cc_\alpha$ such that
\begin{equation}\eqlabel{2.1.1}
\Delta_{\alpha,\beta}(\gamma_{\alpha\beta}(c))=\gamma_\alpha(c_{(1)})\ot_A 
\gamma_\beta(c_{(2)}),
\end{equation}
for all $c\in \Cc$.
\end{definition}

From \equref{1.1.2} and \equref{2.1.1}, it follows that
$$(\varepsilon\circ \gamma_\alpha^{-1})(\gamma_{\alpha\beta}(c)_{(1,\alpha)})\gamma_{\alpha\beta}(c)_{(2,\beta)}=
\gamma_\beta(c),$$
for all $c\in \Cc$. This can be restated as follows: for all $c\in \gamma_{\alpha\beta}(\Cc)=
\Cc_{\alpha\beta}$, we have
\begin{equation}\eqlabel{2.1.2}
(\varepsilon\circ \gamma_\alpha^{-1})(c_{(1,\alpha)})c_{(2,\beta)}=(\gamma_\beta\circ \gamma_{\alpha\beta}^{-1})(c).
\end{equation}
In a similar way, we obtain the formula
\begin{equation}\eqlabel{2.1.3}
c_{(1,\alpha)}(\varepsilon\circ \gamma_\beta^{-1})(c_{(2,\beta)})=(\gamma_\alpha\circ \gamma_{\alpha\beta}^{-1})(c).
\end{equation}
A cofree group coring $\ul{\Cc}$ is defined up to isomorphism by $\Cc_e$. We will write
$\ul{\Cc}=\Cc_e\lan G\ran $.

\begin{theorem}\thlabel{2.2}
If $\ul{\Cc}$ is a cofree group coring, then the categories $\Mm^{\Cc_e}$ and $\Mm^{G,\ul{\Cc}}$
are equivalent.
\end{theorem}

\begin{proof}
We define a functor $F_2:\ \Mm^{\Cc_e}\to \Mm^{G,\ul{\Cc}}$ as follows: $F_2(N)_\alpha=
\nu_\alpha(N)$ is an isomorphic copy of $N$; the coaction maps are
$$\rho_{\alpha,\beta}:\ \nu_{\alpha\beta}(N)\to \nu_\alpha(N)\ot_A \gamma_\beta(\Cc_e),~~
\rho_{\alpha,\beta}(\nu_{\alpha\beta}(n))=\nu_\alpha(n_{[0]})\ot_A \gamma_\beta(n_{[1]}).$$
We also have a functor $G_2:\ \Mm^{G,\ul{\Cc}}\to \Mm^{\Cc_e}$, $G_2(\ul{M})=M_e$,
with coaction $\rho_{e,e}=\rho$. It is then clear that $G_2(F_2(N))=N$, for all
$N\in \Mm^{\Cc_e}$. For $\ul{M}\in \Mm^{G,\ul{\Cc}}$, we have that
$$F_2(G_2(\ul{M}))=(\nu_\alpha(M_e))_{\alpha\in G}.$$
It is clear that the map
$$\varphi_\alpha:\ M_\alpha\to \nu_\alpha(M_e),~~\varphi_\alpha(m)=\nu_\alpha(m_{[0,e]})
\varepsilon(\gamma_\alpha^{-1}(m_{[1,\alpha]}))$$
is right $A$-linear. $\ul{\varphi}=(\varphi_\alpha)_{\alpha\in G}$ is a morphism in
$\Mm^{G,\ul{\Cc}}$, since
\begin{eqnarray*}
&&\hspace*{-2cm}
((\varphi_\alpha\ot_A \Cc_\beta)\circ \rho_{\alpha,\beta})(m)=
\varphi_\alpha(m_{[0,\alpha]})\ot_A m_{[1,\beta]}\\
&=& \nu_\alpha(m_{[0,e]})(\varepsilon\circ \gamma_\alpha^{-1})(m_{[1,\alpha]})
\ot_A m_{[2,\beta]}\\
&\equal{\equref{2.1.2}}&  \nu_\alpha(m_{[0,e]})\ot_A (\gamma_\beta\circ \gamma_{\alpha\beta}^{-1})
(m_{[1,\alpha\beta]})\\
&\equal{\equref{2.1.3}}&
\nu_{\alpha}(m_{[0,e]})\ot_A \gamma_\beta\Bigl(m_{[1,e]}(\varepsilon\circ \gamma_{\alpha\beta}^{-1})
(m_{[2,\alpha\beta]})\Bigr)\\
&=&\nu_{\alpha}(m_{[0,e]})\ot_A \gamma_\beta(m_{[1,e]})
(\varepsilon\circ \gamma_{\alpha\beta}^{-1})(m_{[2,\alpha\beta]})\\
&=& \rho_{\alpha,\beta}(\nu_{\alpha\beta}(m_{[0,e]}))(\varepsilon\circ \gamma_{\alpha\beta}^{-1})(m_{[1,\alpha\beta]})\\
&=& (\rho_{\alpha,\beta}\circ \varphi_{\alpha\beta})(m),
\end{eqnarray*}
for all $m\in M_{\alpha\beta}$. Next we define
$$\psi_\alpha:\ \nu_\alpha(M_e)\to M_\alpha,~~\psi_\alpha(\nu_\alpha(m))=
m_{[0,\alpha]}(\varepsilon\circ \gamma_{\alpha^{-1}}^{-1})(m_{[1,\alpha^{-1}]}).$$
For all $m\in M_\alpha$, we compute
\begin{eqnarray*}
&&\hspace*{-2cm}
(\psi_\alpha\circ\varphi_\alpha)(m)=
\psi_\alpha(\nu_\alpha(m_{[0,e]})
\varepsilon(f_\alpha^{-1}(m_{[1,\alpha]})))\\
&=& m_{[0,\alpha]}(\varepsilon\circ \gamma_{\alpha^{-1}}^{-1})(m_{[1,\alpha^{-1}]})
(\varepsilon\circ \gamma_\alpha^{-1})(m_{[2,\alpha]})\\
&=&m_{[0,\alpha]}(\varepsilon\circ \gamma_{\alpha^{-1}}^{-1})\Bigl(m_{[1,\alpha^{-1}]}
(\varepsilon\circ \gamma_\alpha^{-1})(m_{[2,\alpha]})\Bigr)\\
&\equal{\equref{2.1.3}}& 
m_{[0,\alpha]}(\varepsilon\circ \gamma_{\alpha^{-1}}^{-1})\Bigl(
(\gamma_{\alpha^{-1}}\circ \gamma_e^{-1})(m_{[1,e]})\Bigr)\\
&=& m_{[0,\alpha]}\varepsilon(m_{[1,e]})=m.
\end{eqnarray*}
For all $m\in M_e$, we have
\begin{eqnarray*}
&&\hspace*{-2cm}
(\varphi_\alpha\circ\psi_\alpha)(\nu_\alpha(m))=
\varphi_\alpha\Bigl(m_{[0,\alpha]}(\varepsilon\circ \gamma_{\alpha^{-1}}^{-1})(m_{[1,\alpha^{-1}]})\Bigr)\\
&=& \nu_\alpha(m_{[0,e]})(\varepsilon\circ \gamma_\alpha^{-1})(m_{[1,\alpha]})
(\varepsilon\circ \gamma_{\alpha^{-1}}^{-1})(m_{[2,\alpha^{-1}]})\\
&=& \nu_\alpha(m_{[0,e]})(\varepsilon\circ \gamma_\alpha^{-1})\Bigl(m_{[1,\alpha]}
(\varepsilon\circ \gamma_{\alpha^{-1}}^{-1})(m_{[2,\alpha^{-1}]})\Bigr)\\
&\equal{\equref{2.1.3}}& 
\nu_\alpha(m_{[0,e]}) (\varepsilon\circ \gamma_\alpha^{-1}\circ \gamma_\alpha\circ \gamma_e^{-1})(m_{[1,e]})\\
&=& \nu_\alpha(m_{[0,e]}\varepsilon(m_{[1,e]}))=\nu_\alpha(m).
\end{eqnarray*}
This shows that $\psi_\alpha$ is inverse to $\varphi_\alpha$, and our result follows.
\end{proof}

\section{Graded corings and comodules}\selabel{3}
Let $\Cc$ be an $A$-coring, $\Cc$ is called a $G$-graded $A$-coring if there exists
a direct sum decomposition
$\Cc=\bigoplus_{\alpha\in G} \Cc_\alpha$
as $A$-bimodules such that
$\Delta(\Cc_\alpha)\subset \bigoplus_{\beta\in G} \Cc_{\alpha\beta^{-1}}\ot_A \Cc_\beta$ and
$\varepsilon(\Cc_\alpha)=0~{\rm if}~\alpha\neq e$.\\
If $A$ is a commutative ring, and $ac=ca$, for all $a\in A$ and $c\in \Cc$,
then $\Cc$ is called a $G$-graded coalgebra, cf. \cite{NasTor}.\\
To a $G$-graded $A$-coring, we can associate a $G$-$A$-coring
$\ul{\Cc}=(\Cc_\alpha)_{\alpha\in G}$. The counit is the restriction of $\varepsilon$ to
$\Cc_e$, and $\Delta_{\alpha,\beta}$ is the composition
$$\Cc_{\alpha\beta}\mapright{\Delta}
\bigoplus_{\gamma\in G} \Cc_{\alpha\beta\gamma^{-1}}\ot_A\Cc_\gamma
\mapright{p} \Cc_\alpha\ot_A \Cc_\beta,$$
where $p$ is the obvious projection.\\

Let $(M,\rho)$ be a right $\Cc$-comodule. For each $\alpha\in G$, we consider the map
$$(M\ot_A p_\alpha)\circ\rho:\ M\to M\ot_A\Cc_\alpha,$$
where $p_\alpha:\ \Cc\to \Cc_\alpha$ is the projection. Then $M$ is a right $\ul{\Cc}$-comodule,
and we obtain a functor $\Mm^{\Cc}\to\Mm^{\ul{\Cc}}$.\\
$(M,\rho)$ is called a $G$-graded right $\Cc$-comodule if we have a decomposition
$M=\bigoplus_{\alpha\in G} M_\alpha$
as right $A$-modules, such that
$$\rho(M_\alpha)\subset \bigoplus_{\beta\in G} M_{\alpha\beta^{-1}}\ot_A \Cc_\beta.$$
Now consider the maps
$$\rho_{\alpha,\beta}:\ M_{\alpha\beta}\mapright{\rho} 
\bigoplus_{\gamma\in G} M_{\alpha\beta\gamma^{-1}}\ot_A \Cc_\gamma
\mapright{p} M_\alpha\ot_A \Cc_\beta.$$
Then $\ul{M}=(M_\alpha)_{\alpha\in G}$ is a right $G$-$\ul{\Cc}$-comodule, and
we have a functor from the category of graded $\Cc$-comodules to $\Mm^{G,\ul{\Cc}}$.\\

If $G$ is finite, then there is a one-to-one correspondence between graded corings and
group corings: if $\ul{\Cc}$ is a group coring, then $\bigoplus_{\alpha\in G}\Cc_\alpha$
is a graded coring. In this situation the two functors between ($G$-graded) $\Cc$-comodules and
($G$-)$\ul{\Cc}$-comodules are isomorphisms of categories.

\section{Graded rings and modules}\selabel{4}
Let $A$ be a ring and $R=\bigoplus_{\alpha\in G}R_\alpha$ a $G$-graded ring.
Suppose that we have a ring morphism $i:\ A\to R_e$. Then we call $R$ a $G$-graded
$A$-ring. Every $R_\alpha$ is then an $A$-bimodule, via restriction of scalars, and
the decomposition of $R$ is a decomposition of $A$-bimodules. The category of
$G$-graded right $R$-modules will be denoted by $\Mm^G_R$.\\
Let $\ul{\Cc}$ be a $G$-$A$-coring. For every $\alpha\in G$, 
$R_\alpha={}^*\Cc_{\alpha^{-1}}={}_A\Hom(\Cc_{\alpha^{-1}},A)$
is an $A$-bimodule, with
$$(a\cdot f\cdot b)(c)=f(ca)b,$$
for all $f\in R_\alpha$, $a,b\in A$ and $c\in \Cc_{\alpha^{-1}}$.\\
Take $f_\alpha\in R_\alpha$, $g_\beta\in R_\beta$ and define
$f_\alpha\# g_\beta\in R_{\alpha\beta}$ as the composition
$$\Cc_{(\alpha\beta)^{-1}}\mapright{\Delta_{\beta^{-1},\alpha^{-1}}}
\Cc_{\beta^{-1}}\ot_A \Cc_{\alpha^{-1}}\mapright{\Cc_{\beta^{-1}}\ot_A f_\alpha}
\Cc_{\beta^{-1}}\mapright{g_\beta}A,$$
that is,
$$(f_\alpha\#g_\beta)(c)=g_\beta(c_{(1,\beta^{-1})}f_\alpha(c_{(2,\alpha^{-1})})),$$
for all $c\in \Cc_{(\alpha\beta)^{-1}}$.
This defines maps $m_{\alpha,\beta}:\ R_\alpha\ot_A R_\beta\to R_{\alpha\beta}$,
which make $R=\bigoplus_{\alpha\in G}R_\alpha$ into a $G$-graded $A$-ring.
Let us show that the multiplication is associative: take $h_\gamma\in R_\gamma$
and $c\in \Cc_{(\alpha\beta\gamma)^{-1}}$. We then compute that
\begin{eqnarray*}
&&\hspace*{-2cm}
((f_\alpha\#g_\beta)\#h_\gamma)(c)=
h_\gamma(c_{(1,\gamma^{-1})}(f_\alpha\#g_\beta)(c_{(2,(\alpha\beta)^{-1})}))\\
&=& h_\gamma(c_{(1,\gamma^{-1})}g_\beta(c_{(2,\beta^{-1})}f_\alpha(c_{(3,\alpha^{-1})})))
= (f_\alpha\#(g_\beta\#h_\gamma))(c).
\end{eqnarray*}
$\varepsilon\in R_e$ is a unit for the multiplication; $i:\ A\to R_e$, $i(a)(c)=\varepsilon(c)a$
is a ring homomorphism, since
\begin{eqnarray*}
&&\hspace*{-2cm}
(i(a)\#i(b))(c)=i(b)(c_{(1,e)}i(a)(c_{(2,e)}))=i(b)(c_{(1,e)}\varepsilon(c_{(2,e)})a)\\
&=& i(b)(ca)=\varepsilon(ca)b=\varepsilon(c)ab=i(ab)(c).
\end{eqnarray*}
We conclude that
$R=\bigoplus_{\alpha\in G} R_\alpha$ is a $G$-graded $A$-ring, called the (left) dual (graded) ring of
the group coring $\ul{\Cc}$. We will also write ${}^*\ul{\Cc}=R$.\\
Suppose we are given a morphism $\ul{f}=(f_\alpha)_{\alpha\in G}:\ul{\Cc}\to \ul{\Dd}$ of $G$-$A$-corings, its left dual is defined as the $G$-graded $A$-ring morphism
$${}^*\ul{f}=\bigoplus_{\alpha\in G}{}^* f_{\alpha^{-1}}:{}^*\ul{\Dd}=R'\to {}^*\ul{\Cc}=R,\ {}^*\ul{f}\Big( \sum_{\alpha\in G}g_\alpha\Big)=\sum_{\alpha \in G}g_\alpha\circ f_{\alpha^{-1}}.$$

For every $\alpha\in G$, $R_\alpha^*=\Hom_A(R_\alpha,A)$ is an $A$-bimodule, with
structure maps
$$(a\cdot h\cdot b)(f)=ah(bf),$$
for all $a,b\in A$, $f\in R_\alpha$ and $h\in R_\alpha^*$. We have an $A$-bimodule map
$$\iota_\alpha:\ \Cc_{\alpha^{-1}}\to R_\alpha^*,~~
\iota_\alpha(c)(f)=f(c).$$
If $\Cc_{\alpha^{-1}}$ is finitely generated and projective as a left $A$-module, then
$\iota_\alpha$ is an isomorphism of $A$-bimodules. Since
$$R^*=\Hom_A(R,A)=\Hom_A\Bigl(\bigoplus_{\alpha\in G} R_\alpha,A\Bigr)
\cong \prod_{\alpha\in G} R_\alpha^*\cong \prod_{\alpha\in G} R_{\alpha^{-1}}^*,$$
the $\iota_\alpha$ define an $A$-bimodule map
$$\prod_{\alpha\in G} \iota_{\alpha^{-1}}\cong \iota:\
\prod_{\alpha\in G}\Cc_\alpha\to \prod_{\alpha\in G} R_{\alpha^{-1}}^*\cong R^*.$$
We say that a group $A$-coring $\ul{\Cc}$ is left homogeneously finite if every $\Cc_\alpha$
is finitely generated and projective as a left $A$-module. In this case
$\iota$ is an isomorphism.

\begin{proposition}\prlabel{4.2}
Let $\ul{\Cc}$ be a $G$-$A$-coring, with left dual graded ring $R$. We have a
functor $F_3:\ \Mm^{G,\ul{\Cc}}\to \Mm^G_R$, which is an isomorphism of categories
if $\ul{\Cc}$ is left homogeneously finite.
\end{proposition}

\begin{proof}
Take $\ul{M}=(M_\alpha)_{\alpha\in G}\in \Mm^{G,\ul{\Cc}}$. The maps
$$\psi_{\alpha,\beta}:\ M_\alpha\ot_A R_\beta\to M_{\alpha\beta},~~
\psi_{\alpha,\beta}(m\ot_A f)=m\cdot f=m_{[0,\alpha\beta]}f(m_{[1,\beta^{-1}]}),$$
are well-defined, since
\begin{eqnarray*}
&&\hspace*{-2cm}
\psi_{\alpha,\beta}(ma\ot_A f)=m_{[0,\alpha\beta]}f(m_{[1,\beta^{-1}]}a)\\
&=& m_{[0,\alpha\beta]}(a\cdot f)(m_{[1,\beta^{-1}]})=
\psi_{\alpha,\beta}(m\ot_A a\cdot f).
\end{eqnarray*}
$\psi_{\alpha,\beta}$ is right $A$-linear:
\begin{eqnarray*}
&&\hspace*{-2cm}
\psi_{\alpha,\beta}(m\ot_A f\cdot a)=m_{[0,\alpha\beta]}(f\cdot a)(m_{[1,\beta^{-1}]})\\
&=& m_{[0,\alpha\beta]}f(m_{[1,\beta^{-1}]})a=\psi_{\alpha,\beta}(m\ot_A f)a.
\end{eqnarray*}
We also compute that
$$m\cdot \varepsilon=\psi_{\alpha,e}(m\ot_A \varepsilon)=m_{[0,\alpha]}\varepsilon(m_{[1,e]})
=m;$$
if $g\in R_\gamma$, then we have
\begin{eqnarray*}
&&\hspace*{-15mm}
m\cdot (f\#g)=m_{[0,\alpha\beta\gamma]}(f\#g)(m_{[1,(\beta\gamma)^{-1}]})
= m_{[0,\alpha\beta\gamma]}g(m_{[1,\gamma^{-1}]}f(m_{[2,\beta^{-1}]}))\\
&=& (m_{[0,\alpha\beta]}f(m_{[1,\beta^{-1}]}))\cdot g= (m\cdot f)\cdot g.
\end{eqnarray*}
This shows that $\bigoplus_{\alpha\in G}M_\alpha=F_3(\ul{M})$ is a $G$-graded $R$-module.
For a morphism $\ul{f}:\ \ul{M}\to \ul{N}$ in $\Mm^{G,\ul{\Cc}}$, we define
$F_3(\ul{f})=\bigoplus_{\alpha\in G} f_\alpha$.
\end{proof}

Before we prove the second part of \prref{4.2}, we state and prove two Lemmas.

\begin{lemma}\lelabel{4.3}
Let $A$ be a ring, and $M,P\in {}_A\Mm$, with $M$ finitely generated and projective,
with finite dual basis $f^{(\alpha)}\ot_A m^{(\alpha)}$ (finite sum is implicitly understood).
Then
$$\sum_i f_i\ot_A p_i=\sum_j g_j \ot_A q_j~~{ in }~~{}^*M\ot_AP$$
if and only if
$$\sum_i f_i(m)p_i=\sum_j g_j(m)q_j,$$
for all $m\in M$.
\end{lemma}

\begin{proof}
One direction is obvious. Conversely, we have
\begin{eqnarray*}
&&\hspace*{-2cm}
\sum_i f_i\ot_A p_i=\sum_i f^{(\alpha)}\cdot f_i(m^{(\alpha)})\ot_A p_i\\
&=& \sum_i f^{(\alpha)}\ot_A f_i(m^{(\alpha)})p_i
= \sum_j f^{(\alpha)}\ot_A g_j(m^{(\alpha)})q_j\\
&=& \sum_j f^{(\alpha)}\cdot g_j(m^{(\alpha)})\ot_A q_j= \sum_j g_j \ot_A q_j.
\end{eqnarray*}
\end{proof}

\begin{lemma}\lelabel{4.4}
Let $\ul{\Cc}$ be a left homogeneously finite $G$-$A$-coring. Let
$f^{(\alpha)}\ot_A c^{(\alpha)}\in R_{\alpha^{-1}}\ot_A \Cc_\alpha$ be the
finite dual basis of $\Cc_\alpha$ as a left $A$-module. Then
\begin{equation}\eqlabel{4.4.1}
f^{(\beta\gamma)}\ot_A \Delta_{\beta,\gamma}(c^{(\beta\gamma)})
= f^{(\gamma)}\# f^{(\beta)}\ot_A c^{(\beta)}\ot_A c^{(\gamma)}.
\end{equation}
\end{lemma}

\begin{proof}
For all $c\in \Cc_{\beta\gamma}$, we have
\begin{eqnarray*}
&&\hspace*{-2cm}
(f^{(\gamma)}\# f^{(\beta)})(c) c^{(\beta)}\ot_A c^{(\gamma)}
= f^{(\beta)}(c_{(1,\beta)}f^{(\gamma)}(c_{(2,\gamma)}))c^{(\beta)}\ot_A c^{(\gamma)}\\
&=&c_{(1,\beta)}f^{(\gamma)}(c_{(2,\gamma)})\ot_A c^{(\gamma)}
=c_{(1,\beta)}\ot_A f^{(\gamma)}(c_{(2,\gamma)})c^{(\gamma)}\\
&=& c_{(1,\beta)}\ot_A c_{(2,\gamma)}=\Delta_{\beta,\gamma}(c)\\
&=& \Delta_{\beta,\gamma}( f^{(\beta\gamma)}(c)c^{(\beta\gamma)})
=f^{(\beta\gamma)}(c) \Delta_{\beta,\gamma}( c^{(\beta\gamma)}).
\end{eqnarray*}
\equref{4.4.1} now follows after we apply \leref{4.3}, with $M=\Cc_{\beta \gamma}$ and
$P= \Cc_\beta \ot_A \Cc_\gamma$.
\end{proof}

\begin{proof}[Proof (of the second part of \prref{4.2})]
Assume that every $\Cc_\alpha$ is finitely generated and projective as a left
$A$-module. Let $M$ be a $G$-graded right $R$-module, and consider the maps
$$\rho_{\alpha,\beta}:\ M_{\alpha\beta}\to M_\alpha\ot_A \Cc_\beta,~~
\rho_{\alpha,\beta}(m)=m\cdot f^{(\beta)}\ot_A c^{(\beta)}.$$
These maps make $(M_\alpha)_{\alpha\in G}$ into an object of $\Mm^{G,\ul{\Cc}}$.
We first verify the coassociativity conditions:
\begin{eqnarray*}
&&\hspace*{-15mm}
((M_\alpha\ot_A\Delta_{\beta,\gamma})\circ \rho_{\alpha,\beta\gamma})(m)=
(M_\alpha\ot_A\Delta_{\beta,\gamma})(m\cdot f^{(\beta\gamma)}\ot_A c^{(\beta\gamma)})\\
&=&m\cdot f^{(\beta\gamma)}\ot_A \Delta_{\beta,\gamma}(c^{(\beta\gamma)})
\equal{\equref{4.4.1}}
m\cdot (f^{(\gamma)}\#f^{(\beta)})\ot_A c^{(\beta)}\ot_A c^{(\gamma)}\\
&=&(m\cdot f^{(\gamma)})\cdot f^{(\beta)}\ot_A c^{(\beta)}\ot_A c^{(\gamma)}\\
&=&(\rho_{\alpha,\beta}\ot_A\Cc_\gamma)(m\cdot f^{(\gamma)}\ot_A c^{(\gamma)})
= ((\rho_{\alpha,\beta}\ot_A\Cc_\gamma)\circ \rho_{\alpha\beta,\gamma})(m).
\end{eqnarray*}
The counit property can be verified as follows:
$$((M_\alpha\ot_A\varepsilon)\circ\rho_{\alpha,e})(m)
=m\cdot f^{(e)}\varepsilon(c^{(e)})=m\cdot \varepsilon=m.$$
We have a functor $G_3:\  \Mm^G_R\to \Mm^{G,\ul{\Cc}}$. On the objects, it is defined
by $G_3(M)= \ul{M}=(M_\alpha)_{\alpha\in G}$. For a graded $R$-module map
$f:\ M\to N$, we let $G_3(f)_\alpha:\ M_\alpha\to N_\alpha$ be the restriction of $f$ to
$M_\alpha$.\\
We are done if we can show that $F_3$ and $G_3$ are inverses. First take a graded
$R$-module $M$. Then $(F_3\circ G_3)(M)=\bigoplus_{\alpha\in G} M_\alpha=M$,
with right $R$-action coinciding with the original right $R$-action, since
$$
\psi_{\alpha,\beta}(m\ot_A f)=m_{[0,\alpha\beta]}f(m_{[1,\beta^{-1}]})
=(m\cdot f^{(\beta^{-1})})f(c^{(\beta^{-1})})=m\cdot f.$$
Take $\ul{M}\in \Mm^{G,\ul{\Cc}}$. Then $(G_3\circ F_3)(\ul{M})=G_3(\bigoplus_{\alpha\in G}
M_\alpha)=(M_\alpha)_{\alpha\in G}=\ul{M}$. The coaction maps 
$\tilde{\rho}_{\alpha,\beta}$ on $(G_3\circ F_3)(\ul{M})$
coincide with the coaction maps $\rho_{\alpha,\beta}$ on $\ul{M}$, since
\begin{eqnarray*}
&&\hspace*{-2cm}
\tilde{\rho}_{\alpha,\beta}(m)=m\cdot f^{(\beta)}\ot_A c^{(\beta)}
= m_{[0,\alpha]}f^{(\beta)}(m_{[1,\beta]})\ot_A c^{(\beta)}\\
&=& m_{[0,\alpha]}\ot_A f^{(\beta)}(m_{[1,\beta]})c^{(\beta)}=m_{[0,\alpha]}\ot_A m_{[1,\beta]}
={\rho}_{\alpha,\beta}(m),
\end{eqnarray*}
for all $m\in M_{\alpha\beta}$.
\end{proof}

\begin{proposition}\prlabel{4.5}
Let $\ul{\Cc}$ be a $G$-$A$-coring, with left dual graded ring $R$. We have a
functor $F_4:\ \Mm^{\ul{\Cc}}\to \Mm_R$. $F_4$ is an equivalence of categories if
$\ul{\Cc}$ is left homogeneously finite, and $G$ is a finite group.
\end{proposition}

\begin{proof}
Let $(M,(\rho_\alpha)_{\alpha\in G})\in \Mm^{\ul{\Cc}}$. For each $\alpha\in G$, the map
$$\psi_\alpha:\ M\ot_A R_\alpha \to M,~~\psi_\alpha(m\ot_A f)=m_{[0]}f(m_{[1,\alpha^{-1}]}).$$
is well-defined, since
$$\psi_\alpha(ma\ot_A f)=m_{[0]}f(m_{[1,\alpha^{-1}]}a)=m_{[0]}(a\cdot f)(m_{[1,\alpha^{-1}]})
=\psi_\alpha(m\ot_A a\cdot f).$$
$\psi_\alpha$ is right $A$-linear since
$$\psi_\alpha(m\ot_A f\cdot a)=m_{[0]}(f\cdot a)(m_{[1,\alpha^{-1}]})
=m_{[0]}f(m_{[1,\alpha^{-1}]})a
=\psi_\alpha(m\ot_A f)a.$$
We define a right $R$-action on $M$ as follows:
$$m\cdot f=\sum_{\alpha\in G} \psi_\alpha(m\ot_A f_\alpha),$$
for all $m\in M$ and $f=\sum_{\alpha\in G} f_\alpha\in \bigoplus_{\alpha\in G} R_\alpha$.
This makes $M$ into a right $R$-module, since
$m\cdot \varepsilon=m_{[0]}\varepsilon(m_{[1,e]})=m$ and
\begin{eqnarray*}
&&\hspace*{-2cm}
(m\cdot f)\cdot g= (m_{[0]}f(m_{[1,\alpha^{-1}]}))\cdot g
= m_{[0]}g(m_{[1,\beta^{-1}]}f(m_{[2,\alpha^{-1}]}))\\
&=&m_{[0]}(f\#g)(m_{[1,(\alpha\beta)^{-1}]})=m\cdot (f\# g),
\end{eqnarray*}
for all $f\in R_\alpha$ and $g\in R_\beta$. We let $F_4(M)=M$ with the above action;
if $f:\ M\to N$ is a morphism in $\Mm^{\ul{\Cc}}$, then $f$ is also right $R$-linear, and we
define $F_4(f)=f$.\\
If $G$ is finite and every $\Cc_\alpha$ is finitely generated and projective as a left
$A$-module, then $\Cc=\bigoplus_{\alpha\in G}\Cc_\alpha$ is a ($G$-graded) $A$-coring,
that is finitely generated and projective as a left $A$-module. The left dual of the coring
$\Cc$ is precisely the ring $R$, and
$\Mm^{\ul{\Cc}}\cong \Mm^\Cc\cong \Mm_R$.
\end{proof}

Let $R$ be a $G$-graded ring. It is well-known (see for example \cite[Theorem 2.5.1]{NVO})
that we have a pair of adjoint functors $(F_5,G_5)$ between the categories
$\Mm_R^G$ and $\Mm_R$. $F_5$ is the functor forgetting the $G$-grading; $G_5$
is defined as follows: $G_5(M)=\bigoplus_{\alpha\in G}\mu_\alpha(M)$, with right $R$-action
$$\mu_\alpha(m)r=\mu_{\alpha\beta}(mr),$$
for all $m\in M$ and $r\in R_\beta$.

\begin{proposition}\prlabel{4.6}
Let $\ul{\Cc}$ be a $G$-$A$-coring, with left dual $G$-$A$-ring $R$.
Then we have the following commutative diagram of functors.
$$\xymatrix{
\Mm^{G,\ul{\Cc}}\ar[rr]^{F_3}\ar@<-.5ex>[d]_{F_1}&&\Mm_R^G\ar@<-.5ex>[d]_{F_5}\\
\Mm^{\ul{\Cc}}\ar[rr]_{F_4}\ar@<-.5ex>[u]_{G_1}&&\Mm_R\ar@<-.5ex>[u]_{G_5}
}$$
\end{proposition}

\begin{proposition}\prlabel{4.7}
The left dual of a cofree $G$-$A$-coring $\ul{\Cc}$  is the group ring $R_e[G]$.
\end{proposition}

\begin{proof}
For every $\alpha\in G$, we have $A$-bimodule isomorphisms
$$\gamma_{\alpha^{-1}}:\ \Cc_e\to \Cc_{\alpha^{-1}},$$
$${}^*\gamma_{\alpha^{-1}}:\ {}^*\Cc_{\alpha^{-1}}=R_\alpha\to {}^*\Cc_e=R_e,$$
and
$$\sigma_\alpha= ({}^*\gamma_{\alpha^{-1}})^{-1}={}^*(\gamma_{\alpha^{-1}}^{-1}):\
R_e\to R_\alpha.$$
We then have the following property, for $f\in R_e$ and $c\in \Cc_e$:
$$(\sigma_\alpha(f))\bigl(\gamma_{\alpha^{-1}}(c)\bigr)=
\bigl({}^*(\gamma_{\alpha^{-1}}^{-1})(f)\bigr)\bigl(\gamma_{\alpha^{-1}}(c)\bigr)=
f\Bigr(\gamma_{\alpha^{-1}}^{-1}\bigl(\gamma_{\alpha^{-1}}(c)\bigr)\Bigr)=f(c).$$
Using the formula
$$\Delta_{\beta^{-1},\alpha^{-1}}(\gamma_{(\alpha\beta)^{-1}}(c))=
\gamma_{\beta^{-1}}(c_{(1)})\ot_A \gamma_{\alpha^{-1}}(c_{(2)}),$$
we compute, for all $c\in \Cc_e$ and $f,g\in R_e$ that
\begin{eqnarray*}
&&\hspace*{-2cm}
(\sigma_\alpha(f)\#\sigma_\beta(g))(\gamma_{(\alpha\beta)^{-1}}(c))
=\sigma_\beta(g)\Bigl(\gamma_{\beta^{-1}}(c_{(1)})\sigma_\alpha(f)\bigl(\gamma_{\alpha^{-1}}
(c_{(2)})\bigr)\Bigr)\\
&=&
\sigma_\beta(g)\Bigl(\gamma_{\beta^{-1}}(c_{(1)}) f(c_{(2)})\Bigr)
=
\sigma_\beta(g)\Bigl(\gamma_{\beta^{-1}}(c_{(1)} f(c_{(2)}))\Bigr)\\
&=&g(c_{(1)}f(c_{(2)}))=(f\# g)(c)=\sigma_{\alpha\beta}(f\# g)(\gamma_{(\alpha\beta)^{-1}}(c)).
\end{eqnarray*}
The map $\phi:\ R_e[G]\to R$, $\phi(ru_\alpha)=\sigma_\alpha(r)$ is a bijection. It follows from
the above computations that it preserves the multiplication, so it is an isomorphism of
rings. It is clear that it preserves the grading.
\end{proof}

\section{Galois group corings}\selabel{5}
Let $\ul{\Cc}=(\Cc_\alpha)_{\alpha\in G}$ be a $G$-$A$-coring.
A family $\ul{x}=(x_\alpha)_{\alpha\in G}\in \prod_{\alpha\in G}\Cc_\alpha$
is called grouplike if
$\Delta_{\alpha,\beta}(x_{\alpha\beta})=x_{\alpha}\ot_A x_\beta~~{\rm and}~~
\varepsilon(x_e)=1,$
for all $\alpha,\beta\in G$.

\begin{proposition}\prlabel{5.1}
There is a bijective correspondence between
\begin{itemize}
\item grouplike families of $\ul{\Cc}$;
\item right $\ul{\Cc}$-comodule structures on $A$.
\end{itemize}
\end{proposition}

\begin{proof}
Let $\ul{x}$ be grouplike; the maps
$$\rho_\alpha:\ A\to A\ot_A\Cc_\alpha\cong \Cc_\alpha,~~\rho_\alpha(a)=
1\ot_A x_\alpha a,$$
make $A$ into an object of $\Mm^{\ul{\Cc}}$. Conversely, let $(A,(\rho_\alpha)_{\alpha\in G})
\in \Mm^{\ul{\Cc}}$, and let $x_\alpha=\rho_\alpha(1_A)$. Then $\rho_\alpha(a)=
x_\alpha a$. $(x_\alpha)_{\alpha\in G}$ is grouplike since
\begin{eqnarray*}
&&\hspace*{-2cm}
\Delta_{\alpha,\beta}(x_{\alpha\beta})=((A\ot_A \Delta_{\alpha,\beta})\circ \rho_{\alpha\beta})(1_A)
=
((\rho_\alpha\ot_A\Cc_\beta)\circ \rho_\beta)(1_A)\\
&=&
(\rho_\alpha\ot_A\Cc_\beta)(1_A\ot_A x_\beta)=x_\alpha\ot_A x_\beta;\\
&&\hspace*{-2cm}
\varepsilon(x_e)=((A\ot_A \varepsilon)\circ \rho_e)(1_A)=1_A.
\end{eqnarray*}
\end{proof}

\begin{example}\exlabel{5.2}
Let $\ul{\Cc}$ be a cofree group coring, and take a grouplike element $x\in G(\Cc_e)$.
Then $(\gamma_\alpha(x))_{\alpha\in G}$ is a grouplike family, since
$\Delta_{\alpha,\beta}(\gamma_{\alpha\beta}(x))=\gamma_\alpha(x)\ot_A \gamma_\beta(x)$
and $\varepsilon(\gamma_e(x))=1$.
\end{example}

Let $(\ul{\Cc},\ul{x})$ be a $G$-$A$-coring with a fixed grouplike family. For $M\in \Mm^{\ul{\Cc}}$,
we define
$$M^{{\rm co}\ul{\Cc}}=\{m\in M~|~\rho_\alpha(m)=m\ot_A x_\alpha,~~\forall \alpha\in G\}.$$
Then
$$T=A^{{\rm co}\ul{\Cc}}=\{a\in A~|~ax_\alpha=x_\alpha a,~~\forall \alpha\in G\}$$
is a subring of $A$. If $B\to T$ is a morphism of rings, then we have, for all $m\in M^{{\rm co}\ul{\Cc}}$
and $b\in B$ that
$$\rho_\alpha(mb)=m\ot_Ax_\alpha b=m\ot_A bx_\alpha=mb\ot_A x_\alpha,$$
so $mb\in M^{{\rm co}\ul{\Cc}}$. It follows that $M^{{\rm co}\ul{\Cc}}\in \Mm_B$.

\begin{proposition}\prlabel{5.3}
With notation as above, we have a pair of adjoint functors $(F_6=-\ot_BA,G_6=(-)^{{\rm co}\ul{\Cc}})$
between the categories $\Mm_B$ and $\Mm^{\ul{\Cc}}$.
\end{proposition}

\begin{proof}
Let $N\in \Mm_B$. On $N\ot_B A$, we consider the following coaction maps:
$$\rho_\alpha:\ N\ot_B A\to N\ot_BA\ot_A \Cc_\alpha\cong N\ot_B\Cc_\alpha,~~
\rho_\alpha(n\ot_B a)=n\ot_B x_\alpha a.$$
It is straightforward to show that this makes $N\ot_B A$ into an object of $\Mm^{\ul{\Cc}}$.\\
Take $N\in \Mm_B$ and ${M}\in \Mm^{\ul{\Cc}}$. We have an isomorphism
$$\phi:\ \Hom^{\ul \Cc}(N\ot_B A, M)\to \Hom_B(N, M^{{\rm co}\ul{\Cc}}),$$
given by
$$\phi(f)(n)=f(n\ot_B 1_A)~~;~~\phi^{-1}(g)(n\ot_B a)=g(n)a.$$
\end{proof}

We also have a pair of adjoint functors between $\Mm_B$ and $\Mm^{G,\ul{\Cc}}$. For
$\ul{M}\in \Mm^{G,\ul{\Cc}}$, we define
$$\ul{M}^{{\rm co}\ul{\Cc}}=\{(m_\alpha)_{\alpha\in G}\in \prod_{\alpha\in G} M_\alpha~|~
\rho_{\alpha,\beta}(m_{\alpha\beta})=m_\alpha\ot_A x_\beta,~~\forall \alpha,\beta\in G\}.$$
Then $\ul{M}^{{\rm co}\ul{\Cc}}\in \Mm_B$: if $(m_\alpha)_{\alpha\in G}\in \ul{M}^{{\rm co}\ul{\Cc}}$
and $b\in B$, then
$$\rho_{\alpha,\beta}(m_{\alpha\beta}b)=m_\alpha\ot_A x_\beta b=m_\alpha\ot_A bx_\beta
=m_\alpha b\ot_A x_\beta.$$
In \prref{5.4}, we use the functors $G_1$ and $F_6$ defined in Propositions \ref{pr:1.2} and \ref{pr:5.3}.

\begin{proposition}\prlabel{5.4}
With notation as above, we have a pair of adjoint functors $(F_7=G_1\circ F_6,G_7=(-)^{{\rm co}\ul{\Cc}})$
between the categories $\Mm_B$ and $\Mm^{G,\ul{\Cc}}$.
\end{proposition}

\begin{proof}
Observe first that $F_7(N)= (\mu_\alpha(N\ot_B A))_{\alpha\in G}$, with coaction maps
$\rho_{\alpha,\beta}:\ \mu_{\alpha\beta}(N\ot_B A)\to \mu_\alpha(N\ot_B A)\ot_A\Cc_\beta$
given by the formula
$$\rho_{\alpha,\beta}(\mu_{\alpha\beta}(n\ot_B a))=\mu_\alpha(n\ot_B 1_A)\ot_A x_\beta a.$$
Take $N\in \Mm_B$, $\ul{M}\in \Mm^{G,\ul{\Cc}}$.
We define a map
$$\phi:\ \Hom^{G,\ul{\Cc}}(F_7(N),\ul{M})\to \Hom_B(N,\ul{M}^{{\rm co}\ul{\Cc}})$$
as follows. For a morphism $\ul{f}=(f_\alpha)_{\alpha\in G}$ from $F_7(N)$ to $\ul{M}$
in $\Mm^{G,\ul{\Cc}}$, let
$$\phi(\ul{f})(n)=(f_\alpha(\mu_\alpha(n\ot_B 1_A)))_{\alpha\in G}.$$
Then $\phi(\ul{f})(n)\in \ul{M}^{{\rm co}\ul{\Cc}}$, since
\begin{eqnarray*}
&&\hspace*{-2cm}
\rho_{\alpha,\beta}(f_{\alpha\beta}(\mu_{\alpha\beta}(n\ot_B 1_A)))=
(f_\alpha\ot_A\Cc_\beta)(\rho_{\alpha,\beta}(\mu_{\alpha\beta}(n\ot_B 1_A)))\\
&=&
(f_\alpha\ot_A\Cc_\beta)(\mu_\alpha(n\ot_B 1_A)\ot_A x_\beta)\\
&=&f_\alpha(\mu_\alpha(n\ot_B 1_A))\ot_A x_\beta.
\end{eqnarray*}
We then define
$$\psi:\ \Hom_B(N,\ul{M}^{{\rm co}\ul{\Cc}})\to \Hom^{G,\ul{\Cc}}(F_7(N),\ul{M})$$
as follows: for $g:\ N\to \ul{M}^{{\rm co}\ul{\Cc}}\subset \prod_{\alpha\in G} M_\alpha$ and $n\in N$, we write
$g(n)=(g(n)_\alpha)_{\alpha\in G}$.
Then we put
$$\psi(g)_\alpha(\mu_\alpha(n\ot_B a))=g(n)_\alpha a.$$
Let us show that $\psi(g)$ is a morphism in $\Mm^{G,\ul{\Cc}}$.
\begin{eqnarray*}
&&\hspace*{-2cm}
\bigl((\psi(g)_\alpha \ot_A \Cc_\beta)\circ \rho_{\alpha,\beta}\bigr)(\mu_{\alpha\beta}(n\ot_B a))\\
&=& (\psi(g)_\alpha \ot_A \Cc_\beta) ( \mu_\alpha(n\ot_B 1_A)\ot_A x_\beta a)
= g(n)_\alpha \ot_A x_\beta a\\
&=& \rho_{\alpha,\beta}(g(n)_{\alpha\beta}a)
= (\rho_{\alpha,\beta}\circ \psi(g)_{\alpha\beta})(\mu_{\alpha\beta}(n\ot_B a)).
\end{eqnarray*}
We can easily show that $\psi$ is inverse to $\phi$. First, $(\psi\circ \phi)(\ul{f})=\ul{f}$, since
$
(\psi\circ \phi)(\ul{f})_\alpha(\mu_\alpha(n\ot_B a))
=f_\alpha(\mu_\alpha(n\ot_B 1_A))a
=f_\alpha(\mu_\alpha(n\ot_B 1_A)a)=f_\alpha( \mu_\alpha(n\ot_B a))$.\\
Secondly, $(\phi\circ\psi)(g)=g$, since
$(\phi\circ\psi)(g)(n)=(\psi(g)_\alpha(\mu_\alpha(n\ot_B 1_A)))_{\alpha\in G}=(g(n)_\alpha)_{\alpha\in G}
=g(n)$.
\end{proof}

\begin{remark}
If $G$ is a finite group we can obtain $(F_7,G_7)$ as the composition of the two pairs of adjoint functors $(G_1,F_1)$ and $(F_6,G_6)$ (see Propositions \ref{pr:1.2} and \ref{pr:5.3}): $(F_7=G_1\circ F_6,G_7=G_6\circ F_1)$. Indeed, let us show that, for $\ul{M}\in \Mm^{G,\ul{\Cc}}$,
$$(G_6\circ F_1)(\ul{M})=\Big(\bigoplus_{\alpha \in G}M_\alpha\Big)^{{\rm co}\ul{\Cc}}$$
equals $\ul{M}^{\rm{co}\ul{\Cc}}$:
$m=(m_\alpha)_{\alpha \in G}$ is in $(\bigoplus_{\alpha \in G}M_\alpha )^{{\rm co}\ul{\Cc}}$ if and only if $\rho_\beta(m)=m\ot_A x_\beta$, for all $\beta \in G$, if and only if,
$$\sum_{\alpha\in G}{m_\alpha}_{[0,\alpha \beta^{-1}]} \ot_A {m_\alpha}_{[1,\beta]}=\sum_{\alpha \in G}m_\alpha \ot_A x_\beta =\sum_{\alpha \in G}m_{\alpha\beta^{-1}} \ot_A x_\beta,$$
for all $\beta \in G$;
this is equivalent to
$$\rho_{\alpha \beta^{-1}, \beta}(m_\alpha)=m_{\alpha \beta^{-1}}\ot_A x_\beta,$$
or 
$\rho_{\alpha,\beta}(m_{\alpha\beta})=m_\alpha\ot_A x_\beta$, for all $\alpha,\beta \in G$, i.e. $m\in \ul{M}^{\rm{co}\ul{\Cc}}$.
\end{remark}

Our next goal is to investigate when $(F_7,G_7)$ is a pair of inverse equivalences.
To this end, we will need the unit and counit of this adjunction. We first describe the unit $\eta_7$:
$$\eta_{7,N}:\ N\to (\mu_\alpha(N\ot_B A))_{\alpha\in G}^{{\rm co}\ul{\Cc}},~~
\eta_{7,N}(n)=(\mu_\alpha(n\ot_B 1_A))_{\alpha\in G}.$$
The counit $\varepsilon_7$ is the following:
$$\varepsilon_{7,\ul{M},\alpha}:\  \mu_\alpha(\ul{M}^{{\rm co}\ul{\Cc}}\ot_B A)\to M_\alpha,~~
\varepsilon_{7,\ul{M},\alpha}(\mu_\alpha((m_\beta)_{\beta\in G}\ot_B a))=m_\alpha a.$$
We will proceed as in \cite{C03}. Let $\Dd_e=A\ot_B A$ be the Sweedler canonical coring
associated to the ring morphism $B\to A$. Recall from \cite{Br3,BrW,Sweedler65} that its 
comultiplication and counit are given by the formulas
$$\Delta(a\ot_B b)=(a\ot_B 1_A)\ot_A (1_A\ot_B b)~~;~~\varepsilon(a\ot_B b)=ab.$$
Let $\ul{\Dd}=(A\ot_B A)\lan G\ran $ be the cofree group coring built on the Sweedler canonical coring.

\begin{lemma}\lelabel{5.5}
We have a morphism of $G$-$A$-corings $\ul{\can}:\ \ul{\Dd}\to \ul{\Cc}$ given by
$$\can_\alpha(\mu_\alpha(a\ot_B b))=ax_\alpha b.$$
\end{lemma}

\begin{proof}
\begin{eqnarray*}
&&\hspace*{-15mm}
((\can_\alpha\ot_A\can_\beta)\circ \Delta_{\alpha,\beta})(\mu_{\alpha\beta}(a\ot_B b))\\
&=&(\can_\alpha\ot_A\can_\beta)(\mu_\alpha(a\ot_B 1_A)\ot_A \mu_\beta(1_A\ot_B b))\\
&=&
ax_\alpha\ot_Ax_\beta b=\Delta_{\alpha,\beta}(ax_{\alpha\beta}b)
=(\Delta_{\alpha,\beta}\circ\can_{\alpha\beta})(\mu_{\alpha\beta}(a\ot_B b));\\
&&\hspace*{-15mm}
\varepsilon(\can_e(a\ot_B b))=\varepsilon(ax_e b)=ab=\varepsilon(a\ot_B b).
\end{eqnarray*}
\end{proof}

\begin{proposition}\prlabel{5.7}
With notation as in \prref{5.4} and \leref{5.5}, we have the following properties.
\begin{enumerate}
\item If $F_7$ is fully faithful, then $i:\ B\to T$ is an isomorphism;
\item if $G_7$ is fully faithful, then $\ul{\can}:\ \ul{\Dd}\to \ul{\Cc}$ is an isomorphism.
\end{enumerate}
\end{proposition}

\begin{proof}
1) Let $\ul{A}=G_1(A)=(\mu_\alpha(A))_{\alpha\in G}$, with
$$\rho_{\alpha,\beta}(\mu_{\alpha\beta}(a))=\mu_\alpha(1_A)\ot_A x_\beta a.$$
Then $(\mu_\alpha(a_\alpha))_{\alpha\in G}\in \ul{A}^{{\rm co}\ul{\Cc}}$ if and only if
$$\rho_{\alpha,\beta}(\mu_{\alpha\beta}(a_{\alpha\beta}))=
\mu_\alpha(1_A)\ot_A x_\beta a_{\alpha\beta}=\mu_\alpha(a_{\alpha})\ot_A x_\beta,$$
or
\begin{equation}\eqlabel{5.7.1}
x_\beta a_{\alpha\beta}=a_\alpha x_\beta,
\end{equation}
for all $\alpha,\beta \in G$. We have an injective map
$$f:\ T={A}^{{\rm co}\ul{\Cc}}\to \ul{A}^{{\rm co}\ul{\Cc}},~~f(a)=
(\mu_\alpha(a))_{\alpha\in G}.$$
Indeed, if $a\in {A}^{{\rm co}\ul{\Cc}}$, then $ax_\alpha=x_\alpha a$, for all $\alpha\in G$,
and then \equref{5.7.1} holds. If $F_7$ is fully faithful, then $\eta_7$ is a natural isomorphism.
In particular,
$\eta_{7,B}:\ B\to \ul{A}^{{\rm co}\ul{\Cc}}$ is an isomorphism. We have that
$$\eta_{7,B}(b)=(\mu_\alpha(b1_A))_{\alpha\in G}=((\mu_\alpha \circ i)(b))_{\alpha\in G}
=(f\circ i)(b).$$
From the fact that $\eta_{7,B}$ is surjective, it follows that $f$ is surjective, so $f$ is
an isomorphism. Since $\eta_{7,B}=f\circ i$, it follows that $i$ is an isomorphism.\\

2) $\ul{\Cc}\in \Mm^{G,\ul{\Cc}}$, with coaction maps $\Delta_{\alpha,\beta}$. We have an isomorphism
$$f:\ A\to \ul{\Cc}^{{\rm co}\ul{\Cc}},~~f(a)=(ax_\alpha)_{\alpha\in G}.$$
$f(a)\in \ul{\Cc}^{{\rm co}\ul{\Cc}}$ since
$$\Delta_{\alpha,\beta}(ax_{\alpha\beta})=ax_\alpha\ot_A x_\beta.$$
The inverse $g$ of $f$ is defined as follows:
$$g((c_\alpha)_{\alpha\in G})=\varepsilon(c_e).$$
It is clear that $(g\circ f)(a)=a$.
For $\ul{c}=(c_\alpha)_{\alpha\in G}\in \ul{\Cc}^{{\rm co}\ul{\Cc}}$, we have
$\Delta_{\alpha,\beta}(c_{\alpha\beta})=c_\alpha\ot_A x_\beta$, and, in particular,
$\Delta_{e,\beta}(c_{\beta})=c_e\ot_A x_\beta$. It follows from \equref{1.1.3} that
$c_\beta=\varepsilon(c_e)x_\beta$. Then we compute that
$$(f\circ g)(\ul{c})=f(\varepsilon(c_e))=(\varepsilon(c_e)x_\alpha)_{\alpha\in G}=\ul{c}.$$
If $G_7$ is fully faithful, then $\varepsilon_7$ is a natural isomorphism. In particular,
$$\varepsilon_{7,\ul{\Cc},\alpha}:\ \mu_\alpha(\ul{\Cc}^{{\rm co}\ul{\Cc}}\ot_B A)\to \Cc_\alpha,
~~\varepsilon_{7,\ul{\Cc},\alpha}(\mu_\alpha(\ul{c}\ot_B a))=c_\alpha a,$$
is an isomorphism. Now we compute that $\can_\alpha$ equals the composition
$$\mu_\alpha(A\ot_B A)
\mapright{\mu_\alpha(f\ot_B A)} \mu_\alpha(\ul{\Cc}^{{\rm co}\ul{\Cc}}\ot_B A)
\mapright{\varepsilon_{7,\ul{\Cc},\alpha}}\Cc_\alpha.$$
Indeed,
\begin{eqnarray*}
&&\hspace*{-2cm}
(\varepsilon_{7,\ul{\Cc},\alpha}\circ (\mu_\alpha(f\ot_B A)))(\mu_\alpha(a\ot_B b))\\
&=&\varepsilon_{7,\ul{\Cc},\alpha}(\mu_\alpha((ax_\beta)_{\beta\in G}\ot_B b))
=ax_\alpha b = \can_\alpha(\mu_\alpha(a\ot_B b)).
\end{eqnarray*}
\end{proof}

Recall (see e.g. \cite{Br3,BrW,C03}) that an $A$-coring with fixed grouplike element
$(\Cc_e,x_e)$ is called a Galois coring if the map
$$\can:\ A\ot_{A^{{\rm co}\Cc_e}}A\to \Cc_e,~~\can(a\ot_{A^{{\rm co}\Cc_e}} b)=ax_eb$$
is an isomorphism of corings. \prref{5.7} suggests the following definition.

\begin{definition}\delabel{5.8}
Let $(\ul{\Cc},\ul{x})$ be a $G$-$A$-coring with a fixed grouplike family.
We say that $(\ul{\Cc},\ul{x})$ is Galois if
$$\ul{\can}:\ \ul{\Dd}=(A\ot_{A^{{\rm co}\ul{\Cc}}}A)\langle G\rangle\to \ul{\Cc}$$
is an isomorphism of group corings.
\end{definition}

If $(\ul{\Cc},\ul{x})$ is a $G$-$A$-coring with a fixed grouplike family, then it is
clear that
$A^{{\rm co}\ul{\Cc}}\subset A^{{\rm co}\Cc_e}$. We will now show that this inclusion is an
equality if $\ul{\Cc}$ is cofree.
If $\ul{\Cc}=\Cc_e\lan G\ran$ is a cofree $G$-$A$-coring, and $\ul{x}$ a grouplike family
of $\ul{\Cc}$ such that $x_\alpha=\gamma_\alpha(x_e)$, for all $\alpha\in G$,
then we will say that $(\ul{\Cc},\ul{x})$ is a cofree group coring with a fixed grouplike
family.

\begin{lemma}\lelabel{5.9}
Let $(\ul{\Cc},\ul{x})$ be a cofree group coring with a fixed grouplike
family. Then
$A^{{\rm co}\ul{\Cc}}=A^{{\rm co}\Cc_e}$.
\end{lemma}

\begin{proof}
If $a\in A^{{\rm co}\Cc_e}$, then $ax_e=x_ea$, hence for all $\alpha\in G$, we have that
$$ax_\alpha=a\gamma_\alpha(x_e)=\gamma_\alpha(ax_e)=\gamma_\alpha(x_ea)
=\gamma_\alpha(x_e)a=x_\alpha a,$$
and it follows that $a\in A^{{\rm co}\ul{\Cc}}$.
\end{proof}

\begin{proposition}\prlabel{5.10}
For a  $G$-$A$-coring with a fixed grouplike family $(\ul{\Cc},\ul{x})$, the following
statements are equivalent:
\begin{enumerate}
\item $(\ul{\Cc},\ul{x})$ is a Galois group coring;
\item $(\ul{\Cc},\ul{x})$ is a cofree group coring with a fixed grouplike
family, and $(\Cc_e,x_e)$ is a Galois coring.
\end{enumerate}
\end{proposition}

\begin{proof}
$\ul{1)\Rightarrow 2)}$. $\ul{\Cc}$ is cofree, since $\ul{\can}:\ \ul{\Dd}\to \ul{\Cc}$ is an isomorphism,
and $\ul{\Dd}$ is cofree. The isomorphisms $\gamma_\alpha:\ \Cc_e\to \Cc_\alpha$
are obtained as follows: $\gamma_\alpha=\can_\alpha\circ \mu_\alpha\circ \can_e^{-1}$.
In particular, $\gamma_\alpha(x_e)=\can_\alpha(\mu_\alpha(1\ot_{A^{{\rm co}\ul{\Cc}}} 1))=x_\alpha$,
and it follows from \leref{5.9} that $A^{{\rm co}\ul{\Cc}}=A^{{\rm co}\Cc_e}$. From the fact
that $\ul{\can}$ is an isomorphism, it follows that
$\can_e:\ A\ot_{A^{{\rm co}\Cc_e}}A\to \Cc_e$, $\can_e(a\ot b)=ax_e b$
is an isomorphism.\\
$\ul{2)\Rightarrow 1)}$. It follows from \leref{5.9} that $A^{{\rm co}\ul{\Cc}}=A^{{\rm co}\Cc_e}$.
The maps
$\can_\alpha:\ \mu_\alpha(A\ot_{A^{{\rm co}\ul{\Cc}}}A)\to \Cc_\alpha=\gamma_\alpha(\Cc_e)$
are given by
$$
\can_\alpha(\mu_\alpha(a\ot b))=ax_\alpha b=a\gamma_\alpha(x_e)b
=\gamma_\alpha(ax_eb)=(\gamma_\alpha\circ \can_e)(a\ot b),
$$
so $\can_\alpha=\gamma_\alpha\circ \can_e$ is an isomorphism.
\end{proof}

Let $(\Cc_e,x_e)$ be an $A$-coring with fixed grouplike element, and let $i:\
B\to A^{{\rm co}\Cc_e}$ be a ring morphism. Recall from \cite{Br3} or
\cite[Sec. 1]{C03} that we have a pair of adjoint functors
$(F_8=-\ot_B A,G_8=(-)^{{\rm co}\Cc_e})$ between $\Mm_B$ and $\Mm^{\Cc_e}$.
Then the following statements are equivalent (cf. \cite[Prop. 3.1 and 3.8]{C03}):
\begin{enumerate}
\item $B=A^{{\rm co}\Cc_e}$, $(\Cc_e,x_e)$ is Galois, and $A$ is faithfully flat as a left
$B$-module;
\item $(F_8,G_8)$ is a pair of inverse equivalences, and $A$ is flat as a left
$B$-module.
\end{enumerate}

\begin{lemma}\lelabel{5.11}
Let $(\ul{\Cc},\ul{x})$ be a cofree $G$-$A$-coring with a fixed grouplike family.
Let $i:\ B\to A^{{\rm co}\ul{\Cc}}
\cong A^{{\rm co}\Cc_e}$ be a ring morphism. Then $F_7\cong F_2\circ F_8$
and $G_7\cong G_8\circ G_2$. Here $(F_2,G_2)$ are defined as in \thref{2.2}, and
$(F_7,G_7)$ as in \prref{5.4}.
\end{lemma}

\begin{proof}
For $N\in \Mm_B$, we calculate easily that
$$(F_2\circ F_8)(N)=F_2(N\ot_B A)=(\nu_\alpha(N\ot_B A))_{\alpha\in G},$$
with coaction maps
$$\rho_{\alpha,\beta}(\nu_{\alpha\beta}(n\ot_Ba))=
\nu_\alpha(n\ot_B1_A)\ot_A \gamma_\beta(x_ea)=
\nu_\alpha(n\ot_B1_A)\ot_A x_\beta a.$$
We then see that $(F_2\circ F_8)(N)\cong F_7(N)$. From the uniqueness of the
adjoint functor, it then follows that $G_7\cong G_8\circ G_2$.
\end{proof}

\begin{theorem}\thlabel{5.12}
Let $(\ul{\Cc},\ul{x})$ be a $G$-$A$-coring with a fixed grouplike family,
and $i:\ B\to A^{{\rm co}\ul{\Cc}}$ a ring morphism. Then the following assertions
are equivalent.
\begin{enumerate}
\item $B\cong A^{{\rm co}\ul{\Cc}}$, $(\ul{\Cc},\ul{x})$ is a Galois group coring,
and $A$ is faithfully flat as a left $B$-module;
\item $(F_7,G_7)$ is a pair of inverse equivalences between the categories
$\Mm_B$ and $\Mm^{G,\ul{\Cc}}$ and $A$ is flat as a left $B$-module.
\end{enumerate}
\end{theorem}

\begin{proof}
$\ul{1)\Rightarrow 2)}$.
It follows from \prref{5.10} that $\ul{\Cc}$ is cofree, and $x_\alpha=\gamma_\alpha(x_e)$,
and $(\Cc_e,x_e)$ is a Galois coring.
We deduce from \leref{5.9} that $B\cong A^{{\rm co}\ul{\Cc}}=A^{{\rm co}\Cc_e}$.
It follows from \thref{2.2} that $F_2$ is an equivalence, and from the observations
preceding \leref{5.11} that $F_8$ is an equivalence. Consequently
$F_7\cong F_2\circ F_8$ is an equivalence.\\
$\ul{2)\Rightarrow 1)}$.
It follows from \prref{5.7} that $B\cong A^{{\rm co}\ul{\Cc}}$ and that
$(\ul{\Cc},\ul{x})$ is a Galois group coring. From \prref{5.10}, it follows that
$\ul{\Cc}$ is cofree, $x_\alpha=\gamma_\alpha(x_e)$. Then it follows from
\thref{2.2} that $F_2$ is an equivalence. From \leref{5.11}, it follows that $F_7\cong F_2\circ F_8$
is an equivalence, hence $F_8$ is an equivalence. It follows from the
observations preceding \leref{5.11} that $A$ is faithfully flat as a left $B$-module.
\end{proof}

\section{Graded Morita contexts}\selabel{6}
Let $R$ be a $G$-graded ring, and $M,N\in \Mm_R^G$. A right $R$-linear map
$f:\ M\to N$ is called homogeneous of degree $\sigma$ if $f(M_\alpha)\subset N_{\sigma\alpha}$,
for all $\alpha\in G$. The additive group of all right $R$-module maps $M\to N$ of degree $\sigma$ is denoted
by $\HOM_R(M,N)_\sigma$, and we let
$$\HOM_R(M,N)=\bigoplus_{\sigma\in G} \HOM_R(M,N)_\sigma.$$
Let $S$ and $R$ be $G$-graded rings. A $G$-graded Morita context connecting $S$ and $R$
is a Morita context $(S,R,P,Q,\varphi,\psi)$ with the following additional structure:
$P$ and $Q$ are graded bimodules, and the maps $\varphi:\ P\ot_R Q\to S$ and
$\psi:\ Q\ot_S P\to R$ are homogeneous of degree $e$. Graded Morita contexts have
been studied in \cite{Boisen,CVO,Marcus}.\\
It is well-known (see \cite[Sec. II.4]{Bass}) that we can associate a Morita context to
a module. This construction can be generalized to the graded case as follows.
Let $P$ be a $G$-graded right $R$-module. Then $S=\END_R(P)$ is a $G$-graded ring,
and $Q=\HOM_R(P,R)\in {}_R\Mm_S^G$, with structure
\begin{equation}\eqlabel{6.1.1}
(r\cdot q\cdot s)(p)=rq(s(p)),
\end{equation}
for all $r\in R$, $s\in S$, $q\in Q$ and $p\in P$. The connecting maps are the following
$$\varphi:\ P\ot_R Q\to S,~~\varphi(p\ot_Rq)(p')=pq(p');$$
$$\psi:\ Q\ot_S P\to R,~~\psi(q\ot_S p)=q(p).$$
Straightforward computations then show that $(S,R,P,Q,\varphi,\psi)$ is a graded Morita
context.

\begin{example}\exlabel{6.1}
Let $\MM_e=(S_e,R_e,P_e,Q_e,\varphi_e,\psi_e)$ be a Morita context, and consider
the group rings $S=S_e[G]$ and $R=R_e[G]$. Then $P=P_e[G]=\oplus_{\sigma\in G}P_eu_\sigma
\in {}_S\Mm_R^G$ and $Q=Q_e[G]=\oplus_{\sigma\in G}Q_eu_\sigma
\in {}_R\Mm_S^G$, with
$$(su_\sigma)(pu_\tau)(ru_\rho)=spr u_{\sigma\tau\rho}~{\rm and}~
(ru_\sigma)(qu_\tau)(su_\rho)=rqs u_{\sigma\tau\rho},$$
for all $\sigma,\tau,\rho\in G$, $r\in R_e$, $s\in S_e$, $p\in P_e$, $q\in Q_e$.
We have well-defined maps
$$\varphi:\ P\ot_RQ\to S,~~\varphi(pu_\sigma\ot_R q u_\tau)=\varphi_e(p\ot_{R_e}q)
u_{\sigma\tau};$$
$$\psi:\ Q\ot_SP\to R,~~\psi(qu_\sigma\ot_S p u_\tau)=\psi_e(q\ot_{S_e}p)
u_{\sigma\tau}.$$
Then $\MM_e[G]=(S,R,P,Q,\varphi,\psi)$ is a graded Morita context.\\
Let $P\in \Mm_R^G$ be a graded $R$-module, where $R=R_e[G]$ is a group ring.
By the Structure Theorem for graded modules over strongly graded rings, $P=P_e[G]$.
Let $\MM_e$ be the Morita context associated to the right $R_e$-module $P_e$.
It is then straightforward to verify that $\MM_e[G]$ is the graded Morita context associated
to $P$.
\end{example}

\section{Morita contexts associated to a group coring}\selabel{7}
Let $(\ul{\Cc},\ul{x})$ be a $G$-$A$-coring with a fixed grouplike family.
We have seen in \prref{5.1} that $A\in \Mm^{\ul{\Cc}}$. The map
$$\chi:\ R\to A,~~\chi(f)=\sum_{\alpha\in G} f_\alpha(x_{\alpha^{-1}})$$
is a right grouplike character (see \cite[Sec. 2]{CVW}; the terminology was introduced in
\cite{CG}). This means that $\chi$ satisfies the following properties:
\begin{enumerate}
\item $\chi$ is right $A$-linear;
\item $\chi(\chi(f)\cdot g)=\chi (f\# g)$;
\item $\chi(\varepsilon)=1_A$.
\end{enumerate}
Verification of these properties is left to the reader; using the right grouplike character
$\chi$ or \prref{4.5}, we find that $A\in \Mm_R$, with structure
$a\rightact f=f(x_{\beta^{-1}}a)$,
for $f\in R_\beta$. We have already considered the subring
$T=A^{{\rm co}\ul{\Cc}}$ of $A$. $T$ is a subring of
$$T'=A^R=\{a\in A~|~f_{\alpha^{-1}}(ax_\alpha)=f_{\alpha^{-1}}(x_\alpha a),~
\forall \alpha\in G,\forall f_{\alpha^{-1}}\in R_{\alpha^{-1}}\}.$$
If $\ul{\Cc}$ is left homogeneously finite, then $T=T'$: if $a\in T'$, then we find, using
the same notation as in \leref{4.4}, that
$$ax_\alpha=f^{(\alpha)}(ax_\alpha)c^{(\alpha)}=f^{(\alpha)}(x_\alpha a)c^{(\alpha)}=
x_\alpha a,$$
so $a\in T$.
From \cite[Prop. 2.2]{CVW}, it follows that we have a Morita context $\MM'=
(T',R,A,O',\tau',\mu')$. Recall that
$$O'=R^R=\{q\in R~|~q\#f=q\cdot \chi(f),~\forall f\in R\}.$$
This means that $q=\sum_{\alpha\in G} q_\alpha\in O'$ if and only if
$q_\alpha\# f=q_{\alpha\beta}\cdot f(x_{\beta^{-1}})$,
for all $\alpha,\beta\in G$ and $f\in R_\beta$, or, equivalently, for all
$c\in \Cc_{(\alpha\beta)^{-1}}$,
\begin{eqnarray*}
&&\hspace*{-2cm}
f(c_{(1,\beta^{-1})}q_\alpha(c_{(2,\alpha^{-1})}))=(q_\alpha\# f)(c)\\
&=&(q_{\alpha\beta}\cdot f(x_{\beta^{-1}}))(c)=q_{\alpha\beta}(c)f(x_{\beta^{-1}})=
f(q_{\alpha\beta}(c)x_{\beta^{-1}}).
\end{eqnarray*}
We conclude that
\begin{eqnarray*}
&&\hspace*{-2cm}
O'=\{q\in R~|~f(c_{(1,\beta^{-1})}q_\alpha(c_{(2,\alpha^{-1})}))=
f(q_{\alpha\beta}(c)x_{\beta^{-1}}),\\
&&\hspace*{3cm}\forall \alpha,\beta\in G,~f\in R_\beta,~
c\in \Cc_{(\alpha\beta)^{-1}}\}.
\end{eqnarray*}
The connecting maps are the following:
$$\tau':\ A\ot_R O'\to T',~~\tau'(a\ot_R q)=\sum_{\alpha\in G} q_\alpha(x_{\alpha^{-1}}a);$$
$$\mu':\ O'\ot_{T'} A\to R,~~\mu'(q\ot_{T'}a)=q\cdot a.$$
Now consider
$$O=\{q\in R~|~c_{(1,\beta^{-1})}q_\alpha(c_{(2,\alpha^{-1})})=
q_{\alpha\beta}(c)x_{\beta^{-1}},~\forall \alpha,\beta\in G,~
c\in \Cc_{(\alpha\beta)^{-1}}\}.$$
It is clear that $O\subset O'$, and that $O=O'$ if $\ul{\Cc}$ is left homogeneously finite.

\begin{proposition}\prlabel{7.1}
We have a Morita context $\MM=(T,R,A,O,\tau,\mu)$, with $\tau$ and $\mu$ defined by
$$\tau:\ A\ot_R O\to T,~~\tau(a\ot_R q)=\sum_{\alpha\in G} q_\alpha(x_{\alpha^{-1}}a);$$
$$\mu:\ O\ot_{T} A\to R,~~\mu(q\ot_{T}a)=q\cdot a.$$
If $\ul{\Cc}$ is left homogeneously finite, then the Morita contexts $\MM$ and $\MM'$ are
isomorphic.
\end{proposition}

\begin{proof}
We will show that $O$ is a left ideal of $R$. All the other verifications are straightforward.
Take $\alpha,\beta,\gamma\in G$, $q\in O$ and $f\in R_\gamma$. $f\# q\in O$
since we have, for all $c\in \Cc_{(\gamma\alpha\beta)^{-1}}$:
\begin{eqnarray*}
&&\hspace*{-2cm}
c_{(1,\beta^{-1})}(f\# q_\alpha)(c_{(2,(\gamma\alpha)^{-1})})=
c_{(1,\beta^{-1})} q_\alpha(c_{(2,\alpha^{-1})}f(c_{(3,\gamma^{-1})}))\\
&=& q_{\alpha\beta}(c_{(1,(\alpha\beta)^{-1})}f(c_{(2,\gamma^{-1})}))x_{\beta^{-1}}
=(f\# q_{\alpha\beta})(c)x_{\beta^{-1}}.
\end{eqnarray*}
\end{proof}

\section{Graded Morita contexts associated to a group coring}\selabel{8}
Let $(\ul{\Cc},\ul{x})$ be a $G$-$A$-coring with a fixed grouplike family. Then $A\in \Mm^{\ul{\Cc}}$
(see \prref{5.1}). From \prref{1.2}, it follows that
$G_1(A)=(\mu_\alpha(A))_{\alpha\in G}\in \Mm^{G,\ul{\Cc}}$,
with coaction maps
$$\rho_{\alpha,\beta}:\ \mu_{\alpha\beta}(A)\to \mu_\alpha(A)\ot_A\Cc_\beta,~~
\rho_{\alpha,\beta}(\mu_{\alpha\beta}(a))=\mu_\alpha(1_A)\ot_A x_\beta a.$$
From \prref{4.2}, we then obtain that
$$A\{G\}=F_3G_1(A)=\bigoplus_{\alpha\in G}\mu_\alpha(A)\in \Mm_R^G.$$
The right $R$-action is defined by the following formula, for $f\in R_\beta$:
\begin{equation}\eqlabel{8.1.1}
\mu_\alpha(a)\cdot f=\mu_{\alpha\beta}(f(x_{\beta^{-1}}a)).
\end{equation}
We will compute the graded Morita context associated to the graded right $R$-module
$A\{G\}$. Consider the rings
$$S'=\{\ul{b}=(b_{\alpha})_{\alpha\in G}\in \prod_{\alpha\in G} A~|~
f(b_{\alpha\beta}x_{\beta^{-1}})=f(x_{\beta^{-1}}b_\alpha),~
\forall \alpha,\beta\in G,~f\in R_\beta\},$$
$$S=\{\ul{b}=(b_{\alpha})_{\alpha\in G}\in \prod_{\alpha\in G} A~|~
b_{\alpha\beta}x_{\beta^{-1}}=x_{\beta^{-1}}b_\alpha,~
\forall \alpha,\beta\in G\}.$$
Clearly $S\subset S'$, and $S=S'$ if $\ul{\Cc}$ is left homogeneously finite.
Observe that we have ring monomorphisms
$$i:\ T\to S,~~i(b)=(b)_{\alpha\in G},~~
i':\ T'\to S',~~i'(b)=(b)_{\alpha\in G}.$$
 On $S$ and $S'$, we have the following right $G$-action:
\begin{equation}\eqlabel{8.1.2}
\ul{b}^\sigma=(b_{\sigma\alpha})_{\alpha\in G}.
\end{equation}
Indeed, if $\ul{b}\in S$, then $\ul{b}^\sigma\in S$, since
$b_{\sigma\alpha\beta}x_{\beta^{-1}}=x_{\beta^{-1}}b_{\sigma\alpha}$, for all
$\alpha,\beta\in G$. In the same way, ${S'}^\sigma\subset S'$. 

\begin{lemma}\lelabel{8.1}
With notation as above, we have that $S^G=T$ and $S'^G=T'$.
\end{lemma}

\begin{proof}
Using the monomorphism $i$, we find immediately that $T\subset S^G$.
If $\ul{b}\in S^G$, then for all $\sigma\in G$, we have that
$(b_{\sigma\alpha})_{\alpha\in G}=(b_{\alpha})_{\alpha\in G}$, hence
$b_\sigma=b_e$, and $\ul{b}=i(b_e)$.
\end{proof}

Now we consider the twisted group rings
$G*S=\bigoplus_{\alpha\in G}u_\alpha S$ and $G*S'=\bigoplus_{\alpha\in G}u_\alpha S'$,
with multiplication
$$u_\alpha\ul{b}u_\beta\ul{c}=u_{\alpha\beta}\ul{b}^\beta\ul{c}.$$

\begin{proposition}\prlabel{8.2}
We have a graded ring isomorphism
$\Xi:\ \END_R(A\{G\})\to G*S'$.
\end{proposition}

\begin{proof}
For every $\sigma\in G$, we have an additive bijection
$\Xi_\sigma:\ \END_R(A\{G\})_\sigma\to u_\sigma S'$, $\Xi_\sigma(h)=u_\sigma \ul{b}$,
with
\begin{equation}\eqlabel{8.2.1}
b_\alpha=(\mu_{\sigma\alpha}^{-1}\circ h\circ \mu_\alpha)(1_A).
\end{equation}
$h$ is completely determined by $\ul{b}$, since
\begin{equation}\eqlabel{8.2.2}
h(\mu_\alpha(a))=h(\mu_\alpha(1_A))a=\mu_{\sigma\alpha}(b_\alpha)a
=\mu_{\sigma\alpha}(b_\alpha a),
\end{equation}
for all $a\in A$. Since $h$ is right $R$-linear, we have, for all $\beta\in G$ and $f\in R_\beta$
that
$$h(\mu_\alpha(1_A)\cdot f)=h(\mu_{\alpha\beta}(f(x_{\beta^{-1}}))\equal{\equref{8.2.2}}
\mu_{\sigma\alpha\beta}(f(b_{\alpha\beta}x_{\beta^{-1}}))$$
equals
$$h(\mu_\alpha(1_A))\cdot f=\mu_{\sigma\alpha}(b_{\alpha})\cdot f=
\mu_{\sigma\alpha\beta}(f(x_{\beta^{-1}}b_\alpha)),$$
so it follows that $f(b_{\alpha\beta}x_{\beta^{-1}})=f(x_{\beta^{-1}}b_\alpha)$, for all
$\alpha,\beta\in G$ and $f\in R_{\beta^{-1}}$. This means that $\ul{b}\in S'$.\\
The inverse of $\Xi_\sigma$ is defined as follows: given $\ul{b}\in S'$, $\Xi_\sigma^{-1}(u_\sigma
\ul{b})=h$ is defined by \equref{8.2.2}. The proof is finished if we can show that
$$\Xi=\bigoplus_{\sigma\in G}\Xi_\sigma:\ \END_R(A\{G\})\to G*S'$$
preserves the multiplication and the unit. Take $h\in \END_R(A\{G\})_\sigma$,
$k\in \END_R(A\{G\})_\tau$, and let $\Xi_\sigma(h)=u_\sigma \ul{b}$ and
$\Xi_\tau(k)=u_\tau\ul{c}$. Then $k\circ h\in \END_R(A\{G\})_{\tau\sigma}$, and
$\Xi_{\tau\sigma}(k\circ h)=u_{\tau\sigma}\ul{d}$, with
\begin{eqnarray*}
d_\alpha&=& (\mu_{\tau\sigma\alpha}^{-1}\circ k\circ h\circ \mu_\alpha)(1_A)
= (\mu_{\tau\sigma\alpha}^{-1}\circ k\circ\mu_{\sigma\alpha} 
\circ\mu_{\sigma\alpha}^{-1}\circ h\circ \mu_\alpha)(1_A)\\
&=& (\mu_{\tau\sigma\alpha}^{-1}\circ k\circ\mu_{\sigma\alpha} )(b_\alpha)
= (\mu_{\tau\sigma\alpha}^{-1}\circ k\circ\mu_{\sigma\alpha} )(1_A)b_\alpha=c_{\sigma\alpha}b_\alpha.
\end{eqnarray*}
This proves that $\ul{d}=\ul{c}^\sigma\ul{b}$, and $\Xi_{\tau\sigma}(k\circ h)=
u_{\tau\sigma}\ul{c}^\sigma\ul{b}=(u_\tau\ul{c})(u_\sigma\ul{b}).$
Finally, $\Xi_e(A\{G\})=u_e\ul{b}$, with $b_\alpha=(\mu_\alpha^{-1}\circ A\{G\}\circ \mu_\alpha)
(1_A)=1_A$. 
\end{proof}

Our next aim is to describe $\HOM_R(A\{G\},R)$. Consider
\begin{eqnarray*}
&&\hspace*{-15mm}
Q=\{\ul{q}=(q_{\alpha})_{\alpha\in G}\in \prod_{\alpha\in G} R_\alpha~|~
c_{(1,\beta^{-1})}q_\alpha(c_{(2,\alpha^{-1})})=q_{\alpha\beta}(c)x_{\beta^{-1}},\\
&&\hspace*{2cm}\forall \alpha,\beta\in G,~c\in \Cc_{(\alpha\beta)^{-1}}\};\\
&&\hspace*{-15mm}Q'=\{\ul{q}=(q_{\alpha})_{\alpha\in G}\in \prod_{\alpha\in G} R_\alpha~|~
f(c_{(1,\beta^{-1})}q_\alpha(c_{(2,\alpha^{-1})}))=f(q_{\alpha\beta}(c)x_{\beta^{-1}}),\\
&&\hspace*{2cm}\forall \alpha,\beta\in G,~c\in \Cc_{(\alpha\beta)^{-1}},~f\in R_\beta\}.
\end{eqnarray*}
It is clear that $Q\subset Q'$ and that $Q=Q'$ if $\ul{\Cc}$ is left homogeneously finite.

\begin{lemma}\lelabel{8.3}
If $f\in R_\gamma$ and $\ul{q}\in Q$ (resp. $Q'$), then $f\cdot \ul{q}=
(f\# q_{\gamma^{-1}\alpha})_{\alpha\in G}\in Q$ (resp. $Q'$).
\end{lemma}

\begin{proof}
We will prove the first statement; the proof of the second one is similar. For
all $c\in \Cc_{(\alpha\beta)^{-1}}$, we have
\begin{eqnarray*}
&&\hspace*{-2cm}
c_{(1,\beta^{-1})}(f\# q_{\gamma^{-1}\alpha})(c_{(2,\alpha^{-1})})=
c_{(1,\beta^{-1})} q_{\gamma^{-1}\alpha}(c_{(2,\alpha^{-1}\gamma)}f(c_{(3,\gamma^{-1})}))\\
&=&q_{\gamma^{-1}\alpha\beta}(c_{(1,\beta^{-1}\alpha^{-1}\gamma)}}f(c_{(2,\gamma^{-1})}))x_{\beta^{-1}}
=(f\# q_{\gamma^{-1}\alpha\beta})(c)x_{\beta^{-1}.
\end{eqnarray*}
\end{proof}

\begin{lemma}\lelabel{8.4}
If $\ul{q}\in Q$ (resp. $Q'$) and $\ul{b}\in S$ (resp. $S'$), then
$\ul{q}\cdot \ul{b}=(q_\alpha\cdot b_\alpha)_{\alpha\in G}\in Q$ (resp. $Q'$).
\end{lemma}

\begin{proof}
The first statement is easily verified by the following computations:
\begin{eqnarray*}
&&\hspace*{-2cm}
c_{(1,\beta^{-1})}(q_\alpha\cdot b_\alpha)(c_{(2,\alpha^{-1})})=
c_{(1,\beta^{-1})}q_\alpha(c_{(2,\alpha^{-1})})b_\alpha= q_{\alpha\beta}(c)x_{\beta^{-1}}b_\alpha\\
&=& q_{\alpha\beta}(c)b_{\alpha\beta}x_{\beta^{-1}}=(q_{\alpha\beta}\cdot b_{\alpha\beta})
(c)x_{\beta^{-1}}.
\end{eqnarray*}
\end{proof}

\begin{lemma}\lelabel{8.5}
$$QG=\bigoplus_{\alpha\in G}\omega_{\alpha}(Q)\in {}_R\Mm_{G*S}^G~~{\rm and}~~
Q'G=\bigoplus_{\alpha\in G}\omega_{\alpha}(Q')\in {}_R\Mm_{G*S'}^G$$
with bimodule structures defined as follows, for all $f\in R_\beta$, $\ul{q}\in Q$
(resp. $Q'$) and $\ul{b}\in S$ (resp. $S'$):
\begin{equation}\eqlabel{8.5.1}
f\cdot \omega_\alpha(\ul{q})=\omega_{\beta\alpha}(f\cdot \ul{q});
\end{equation}
\begin{equation}\eqlabel{8.5.2}
\omega_\alpha(\ul{q})\cdot u_\tau\ul{b}=\omega_{\alpha\tau}(\ul{q}\cdot {\ul{b}}^{ (\alpha\tau)^{-1}}).
\end{equation}
\end{lemma}

\begin{proposition}\prlabel{8.6}
We have an isomorphism of graded bimodules
$$\Psi:\ \HOM_R(A\{G\},R)\to Q'G.$$
\end{proposition}

\begin{proof}
We have an additive bijection 
$$\Psi_\sigma:\ \HOM_R(A\{G\},R)_\sigma\to \omega_\sigma(Q'),~~\Psi_\sigma(\varphi)=\omega_\sigma(\ul{q}),$$
with $q_\alpha=\varphi(\mu_{\sigma^{-1}\alpha}(1_A))$. $\ul{q}$ determines $\varphi$
completely since
\begin{equation}\eqlabel{8.6.1}
\varphi(\mu_\alpha(a))=\varphi(\mu_\alpha(1_A))\cdot a=q_{\sigma\alpha}\cdot a.
\end{equation}
Take $\beta\in G$ and $f\in R_\beta$. Since $\varphi$ is right $R$-linear, we have that
$$\varphi(\mu_{\sigma^{-1}\alpha}(1_A)\cdot f)=\varphi(\mu_{\sigma^{-1}\alpha\beta}
f(x_{\beta^{-1}}))\equal{\equref{8.6.1}}q_{\alpha\beta}\cdot f(x_{\beta^{-1}})$$
equals
$$\varphi(\mu_{\sigma^{-1}\alpha}(1_A))\# f=q_\alpha\# f.$$
We then have, for all $c\in \Cc_{(\alpha\beta)^{-1}}$ that
$$f(q_{\alpha\beta}(c)x_{\beta^{-1}})=q_{\alpha\beta}(c)f(x_{\beta^{-1}})=
(q_\alpha\# f)(c)=f(c_{(1,\beta^{-1})}q_\alpha (c_{(2,\alpha^{-1})})),$$
and it follows that $\ul{q}\in Q'$. For $\ul{q}\in Q'$, $\Psi_\sigma^{-1}(\omega_\sigma(\ul{q}))
=\varphi$ is defined by \equref{8.6.1}. We now have an additive bijection
$$\Psi=\bigoplus_{\sigma\in G}\Psi_\sigma:\ \HOM_R(A\{G\},R)\to Q'G.$$
$\Psi$ is left $R$-linear: take $f\in R_\beta$, $\varphi\in  \HOM_R(A\{G\},R)_\sigma$, and
put $\Psi_\sigma(\varphi)=\omega_\sigma(\ul{q})$. Then
\begin{eqnarray*}
&&\hspace*{-2cm}
\Psi_{\beta\sigma}(f\cdot \varphi)=\omega_{\beta\sigma}\Bigl((f\cdot \varphi)
(\mu_{\sigma^{-1}\beta^{-1}\alpha}(1_A)\Bigr)_{\alpha\in G}\\
&=& \omega_{\beta\sigma}\Bigl(f\#q_{\beta^{-1}\alpha}\Bigr)_{\alpha\in G}
=f\cdot \omega_\sigma(\ul{q})=f\cdot \Psi_\sigma(\varphi).
\end{eqnarray*}
Finally, $\Psi$ transports the right $\END_R(A\{G\})$-action on $\HOM_R(A\{G\},R)$ to
the right $G*S'$-action on $Q'G$: take $h\in \END_R(A\{G\})_\tau$, and write
$\Xi_\tau(h)=u_\tau\ul{b}$. Then
$$(\varphi\circ h)(\mu_{\tau^{-1}\sigma^{-1}\alpha}(1_A))=
\varphi(\mu_{\sigma^{-1}\alpha}(b_{\tau^{-1}\sigma^{-1}\alpha}))=q_\alpha\cdot
b_{\tau^{-1}\sigma^{-1}\alpha},$$
hence
$$\Psi_{\sigma\tau}(\varphi\circ h)=\omega_{\sigma\tau}(\ul{q}\cdot {\ul{b}}^{(\sigma\tau)^{-1}})
=\omega_\sigma(\ul{q})\cdot u_\tau\ul{b}=\Psi_\sigma(\varphi)\cdot \Xi_\tau(h).$$
\end{proof}

\begin{theorem}\thlabel{8.7}
Consider the graded Morita context $(\END_R(A\{G\}),R,A\{G\},$ $\HOM_R(A\{G\},R),\phi,\psi)$
associated to the graded $R$-module $A\{G\}$. Using the isomorphisms $\Xi$ and $\Psi$
from Propositions \ref{pr:8.2} and \ref{pr:8.6}, we find an isomorphic graded Morita context
$\GG\MM'=(G*S',R,A\{G\},Q'G,\omega',\nu')$, with connecting maps $\omega'$ and $\nu'$
given by the formulas
\begin{eqnarray*}
&&\hspace*{-2cm}
\omega':\ A\{G\}\ot_R Q'G\to G*S',\\
&&\omega'(\mu_\alpha(a)\ot_R \omega_\sigma(\ul{q}))=u_{\alpha\sigma}
\Bigl(q_{\sigma\beta}(x_{(\sigma\beta)^{-1}}a)\Bigr)_{\beta\in G};\\
&&\hspace*{-2cm}
\nu':\ Q'G\ot_{G*S'}A\{G\}\to R,\\
&&\nu'(\omega_\sigma\ul{q}\ot \mu_\alpha(a))=q_{\sigma\alpha}\cdot a.
\end{eqnarray*}
\end{theorem}

\begin{proof}
We have to show that the following two diagrams commute
\begin{equation}\eqlabel{8.7.1}
\diagram
A\{G\}\ot_R \HOM_R(A\{G\},R)\rrto^(.6){\phi}\dto_{A\{G\}\ot_R\Psi}&&\END_R(A\{G\})\dto^{\Xi}\\
A\{G\}\ot_R Q'G\rrto^{\omega'}&&G*S'
\enddiagram
\end{equation}
\begin{equation}\eqlabel{8.7.2}
\diagram
\HOM_R(A\{G\},R)\ot_{\END_R(A\{G\})} A\{G\}\rrto^(.75){\psi}\dto_{\Psi\ot A\{G\}}&&
R\dto^{=}\\
Q'G\ot_{G*S'}A\{G\}\rrto^(.6){\nu'}&&R
\enddiagram
\end{equation}
For all $\alpha,\beta,\sigma\in G$, $a,b\in A$, $\varphi\in \HOM_R(A\{G\},R)_\sigma$,
we have
$$
\phi(\mu_\alpha(a)\ot_R\varphi)(\mu_\beta(b))
=\mu_\alpha(a)\cdot \varphi(\mu_\beta(b))
\equal{\equref{8.1.1}}\mu_{\alpha\sigma\beta}(\varphi(\mu_\beta(b))
(x_{(\sigma\beta)^{-1}}a)),$$
so $(\Xi\circ\phi)(\mu_\alpha(a)\ot_R\varphi)=u_{\alpha\sigma}\ul{b}$, with
$$b_\beta=\mu_{\alpha\sigma\beta}^{-1}
(\phi(\mu_\alpha(a)\ot_R\varphi)(\mu_\beta(1_A))=\varphi(\mu_\beta(1_A))
(x_{(\sigma\beta)^{-1}}a).$$
Now $\mu_\alpha(a)\ot_R \Psi(\varphi)=\mu_\alpha(a)\ot_R \omega_\sigma(\ul{q})$, with
$q_\beta=\varphi(\mu_{\sigma^{-1}\beta}(1_A))$, so
$\omega'(\mu_\alpha(a)\ot_R \Psi(\varphi))=u_{\alpha\sigma}\ul{c}$, with
$$c_\beta=q_{\sigma\beta}(x_{(\sigma\beta)^{-1}}a)=
\varphi(\mu_\beta(1_A))(x_{(\sigma\beta)^{-1}}a)=b_\beta.$$
This proves the commutativity of \equref{8.7.1}. \equref{8.7.2}
also commutes: let $q_\alpha=\varphi(\mu_{\sigma^{-1}\alpha}
(1_A))$. Then
\begin{eqnarray*}
&&\hspace*{-2cm}
(\nu'\circ (\Psi\ot A\{G\}))(\varphi\ot \mu_\alpha(a))=
\nu'(\omega_\sigma(\ul{q})\ot \mu_\alpha(a))=q_{\sigma\alpha}\cdot a\\
&=& \varphi(\mu_\alpha(1_A))\cdot a=\varphi(\mu_\alpha(a))=\psi(\varphi\ot \mu_\alpha(a)).
\end{eqnarray*}
\end{proof}

\begin{theorem}\thlabel{8.8}
Let $(\ul{\Cc},\ul{x})$ be a $G$-$A$-coring with a fixed grouplike family. We have
a second graded Morita context 
$\GG\MM=(G*S,R,A\{G\},QG,\omega,\nu)$, with connecting maps $\omega$ and $\nu$
given by the formulas
\begin{eqnarray*}
&&\hspace*{-2cm}
\omega:\ A\{G\}\ot_R QG\to G*S,\\
&&\omega(\mu_\alpha(a)\ot_R \omega_\sigma(\ul{q}))=u_{\alpha\sigma}
\Bigl(q_{\sigma\beta}(x_{(\sigma\beta)^{-1}}a)\Bigr)_{\beta\in G};\\
&&\hspace*{-2cm}
\nu:\ QG\ot_{G*S}A\{G\}\to R,\\
&&\nu(\omega_\sigma\ul{q}\ot \mu_\alpha(a))=q_{\sigma\alpha}\cdot a.
\end{eqnarray*}
If $\ul{\Cc}$ is left homogeneously finite, then the Morita contexts $\GG\MM'$ and $\GG\MM$
are isomorphic.
\end{theorem}

Let $({\Cc}_e,x_e)$ be an $A$-coring with a fixed grouplike element.
Recall also from \cite{CVW} that we have two Morita contexts
$$\MM_e=(T_e,R_e,A,Q_e,\varphi_e,\psi_e)~~{\rm and}~~
\MM'_e=(T'_e,R_e,A,Q'_e,\varphi'_e,\psi'_e),$$
where
(see \cite{CVW})
$$Q_e=\{q\in R_e={}^*\Cc_e~|~c_{(1)}q(c_{(2)})=q(c)x_e,~\forall c\in \Cc_e\};$$
$$Q'_e=\{q\in R_e={}^*\Cc_e~|~f(c_{(1)}q(c_{(2)}))=f(q(c)x_e),~\forall c\in \Cc_e~{\rm and}~
f\in R_e\};$$
$$T_e=A^{{\rm co}\Cc_e}=\{a\in A~|~ax_e=x_ea\};$$
$$T'_e=A^{R_e}=\{a\in A~|~f(ax_e)=f(x_ea),~\forall f\in R_e\};$$
$$\varphi_e:\ A\ot_{R_e} Q_e\to T_e,~~\varphi_e(a\ot_{R_e}q)=q(x_ea);$$
$$\psi_e:\ Q_e\ot_{T_e} A\to R_e,~~\psi_e(q\ot_{T_e} a)=q\cdot a.$$
$\varphi'_e$ and $\psi'_e$ are defined in a similar way. There is a morphism from
$\MM_e$ to $\MM'_e$, which is an isomorphism if $\Cc_e$ is finitely generated and projective
as a left $A$-module.

\begin{proposition}\prlabel{8.9}
Let $(\ul{\Cc},\ul{x})$ be a cofree group coring with a fixed grouplike
family. Then we have an isomorphism of $G$-graded $R$-modules
$$\vartheta:\ A\{G\}\to A[G],~~\vartheta(\mu_\alpha(a))=au_\alpha.$$
Consequently the graded Morita context $\GG\MM'$ is isomorphic to
$\MM'_e[G]$.
\end{proposition}

\begin{proof}
We have to show that $\vartheta$ is right $R$-linear. Take $f\in R_e={}^*{\Cc_e}$,
$\sigma_\beta(f)\cong f u_\beta\in R_\beta\cong R_e u_\beta$ 
(see \prref{4.7}). Then
\begin{eqnarray*}
&&\hspace*{-15mm}
\mu_\alpha(a)\cdot \sigma_\beta(f)\equal{\equref{8.1.1}}
\mu_{\alpha\beta}(\sigma_\beta(f)(x_{\beta^{-1}}a))
= \mu_{\alpha\beta}((f\circ \gamma_{\beta^{-1}}^{-1})(x_{\beta^{-1}}a))\\
&=& \mu_{\alpha\beta}(f((\gamma_{\beta^{-1}}^{-1}\circ \gamma_{\beta^{-1}})(x_e)a))
= \mu_{\alpha\beta}(f(x_ea)),
\end{eqnarray*}
hence
$$\vartheta(\mu_\alpha(a)\cdot \sigma_\beta(f))=
f(x_ea)u_{\alpha\beta}=(a\rightact  f)u_{\alpha\beta}=(au_\alpha)\cdot (fu_\beta)
=\vartheta(\mu_\alpha(a))\cdot \sigma_\beta(f).$$
The second statement then follows from \exref{6.1}.
\end{proof}

Our next aim is to show that the graded Morita contexts $\GG\MM$ and $\MM_e[G]$
are also isomorphic.

\begin{proposition}\prlabel{8.10}
Let $(\ul{\Cc},\ul{x})$ be a cofree group coring with a fixed grouplike
family. Then $i:\ T\to S$ and $i':\ T'\to S'$ are isomorphisms, and $\END_R(A\{G\})\cong
G*S'$ (resp. $G*S$)
is isomorphic as a graded ring to the group ring $T'[G]$ (resp. $T[G]$).
\end{proposition}

\begin{proof}
It suffices to show that $i$ and $i'$ are surjective. First take $\ul{b}\in S$. Then we
have for $\alpha,\beta\in G$ that
$$\gamma_{\beta^{-1}}(b_{\alpha\beta}x_e)=b_{\alpha\beta}\gamma_{\beta^{-1}}(x_e)
= \gamma_{\beta^{-1}}(x_e) b_\alpha= \gamma_{\beta^{-1}}(x_eb_\alpha).$$
Applying $\varepsilon\circ \gamma_{\beta^{-1}}^{-1}$ to both sides, we find that
$b_{\alpha\beta}=b_{\alpha}$, hence $b_\alpha=b_e$, for all $\alpha\in G$,
and $\ul{b}=i(b_e)$. In a similar way, we find for $\ul{b}\in S'$,
$\alpha,\beta\in G$ and $f\in R_\beta={}^*\Cc_{\beta^{-1}}$ that
$f(\gamma_{\beta^{-1}}(b_{\alpha\beta}x_e))=f(\gamma_{\beta^{-1}}(x_eb_\alpha))$.
Taking $f=\varepsilon\circ \gamma_{\beta^{-1}}^{-1}$, we find that $\ul{b}=i'(b_e)$.
\end{proof}



We have that $T\subset T_e$ and $T'\subset T'_e$. As follows from \leref{5.9}, these inclusions are equalities if
$(\ul{\Cc},\ul{x})$ is a cofree group coring with a fixed grouplike family.

\begin{proposition}\prlabel{8.13}
Let $(\ul{\Cc},\ul{x})$ be a cofree group coring with a fixed grouplike
family. Then $Q\cong Q_e$ and $Q'\cong Q'_e$. Consequently
$\Hom_R(A\{G\},R)\cong Q'_e[G]$.
\end{proposition}

\begin{proof}
We recall from \prref{4.7} that
$R=\bigoplus_{\alpha\in G} \sigma_\alpha(R_e)$, with $\sigma_\alpha(f)=
f\circ \gamma^{-1}_{\alpha^{-1}}$, for all $f\in R_e$. Now take $\ul{q}=
(\sigma_\alpha(q_\alpha))_{\alpha\in G}\in \prod_{\alpha\in G}R_\alpha$. Then for all
$\alpha,\beta\in G$ and $c\in \Cc_e$
$$\gamma_{\beta^{-1}}(c_{(1)})(q_\alpha\circ \gamma^{-1}_{\alpha^{-1}}\circ 
\gamma_{\alpha^{-1}})(c_{(2)})= (q_{\alpha\beta}\circ \gamma^{-1}_{(\alpha\beta)^{-1}}
\circ \gamma_{(\alpha\beta)^{-1}})(c)\gamma_{\beta^{-1}}(x_e),$$
or
$$\gamma_{\beta^{-1}}(c_{(1)})q_\alpha(c_{(2)})=
q_{\alpha\beta}(c) \gamma_{\beta^{-1}}(x_e),$$
or
$$c_{(1)}q_\alpha(c_{(2)})=q_{\alpha\beta}(c)x_e.$$
Taking $\alpha=\beta=e$, we find that $q_e\in Q_e$. Applying $\varepsilon$ to both sides,
we find that $q_\alpha(c)=q_{\alpha\beta}(c)$, and $q_\alpha=q_e$, for all $\alpha\in G$.
These arguments show that the map
$$j:\ Q_e\to Q,~~j(q)=(\sigma_\alpha(q))_{\alpha\in G}$$
is a well-defined isomorphism. In a similar way, we prove that $Q'\cong Q'_e$.
\end{proof}

\begin{theorem}\thlabel{8.12}
Let $(\ul{\Cc},\ul{x})$ be a cofree group coring with a fixed grouplike
family. Then the graded Morita contexts $\GG\MM$ and $\MM_e[G]$ are isomorphic.
\end{theorem}

\begin{proof}
Let $\Theta:\ T[G]\to G*S$ be the isomorphism from \prref{8.10}. We will first show that
the diagram
$$\xymatrix{
A\{G\}\ot_R QG \ar[rr]^{\omega}\ar[d]_{\vartheta\ot j^{-1}G}&&
G*S\\
A[G]\ot_{R_e[G]}Q_e[G]\ar[rr]^{\varphi}&&T_e[G]=T[G]\ar[u]_{\Theta}
}$$
commutes. For $\alpha,\sigma\in G$, $a\in A$ and $\ul{q}\in Q$, we have
\begin{eqnarray*}
&&\hspace*{-15mm}
(\Theta\circ \varphi\circ (\vartheta\ot j^{-1}G))(\mu_\alpha(a)\ot \omega_\sigma(\ul{q}))\\
&=& (\Theta\circ \varphi)(au_\alpha\ot q_eu_\sigma)
= \Theta(q_e(x_ea)u_{\alpha\sigma})\\
&=& u_{\alpha\sigma}\bigl( q_e(x_ea)\bigr)_{\beta\in G}
= u_{\alpha\sigma}\bigl( (q_e\circ \gamma^{-1}_{(\sigma\beta)^{-1}}\circ \gamma_{(\sigma\beta)^{-1}})(x_ea)\bigr)_{\beta\in G}\\
&=&u_{\alpha\sigma}\bigl(q_{\sigma\beta}(x_{(\sigma\beta)^{-1}}a)\bigr)_{\beta\in G}
= \omega(\mu_\alpha(a)\ot \omega_\sigma(\ul{q})).
\end{eqnarray*}
Let $\phi:\ R_e[G]\to R$, $\phi(fu_\alpha)=f\circ \gamma_{\alpha^{-1}}^{-1}$ be the
isomorphism from \prref{4.7}. We will show that the diagram
$$\xymatrix{
QG\ot_{G*S} A\{G\}\ar[rr]^{\nu}\ar[d]_{j^{-1}G\ot \vartheta}&&R\\
Q_e[G]\ot_{T_e[G]}A[G]\ar[rr]^{\psi}&&R_e[G]\ar[u]_{\phi}
}$$
commutes. Take $\alpha,\sigma\in G$, $\ul{q}\in Q$ and $a\in A$. We have seen in
\prref{8.13} that $q_\alpha=q_e\circ \gamma_{\alpha^{-1}}^{-1}$, or
\begin{equation}\eqlabel{8.12.1}
q_e=q_\alpha\circ \gamma_{\alpha^{-1}}.
\end{equation}
We have to show that
\begin{eqnarray*}
&&\hspace*{-2cm}
(\phi\circ\psi\circ (j^{-1}G\ot \vartheta))(\omega_\sigma(\ul{q})\ot \mu_\alpha(a))
= (\phi\circ\psi)(q_eu_\sigma \ot au_\alpha)\\
&=&\phi((q_e\cdot a)u_{\sigma\alpha})=q_e\cdot a\circ \gamma_{(\sigma\alpha)^{-1}}^{-1}
\end{eqnarray*}
and
$$\nu(\omega_\sigma(\ul{q})\ot \mu_\alpha(a))=q_{\sigma\alpha}\cdot a$$
are equal in $R_{\sigma\alpha}={}^*{\Cc}_{(\sigma\alpha)^{-1}}$. For
$\gamma_{(\sigma\alpha)^{-1}}(c)\in \Cc_{(\sigma\alpha)^{-1}}$, we compute that
\begin{eqnarray*}
&&\hspace*{-2cm}
(q_{\sigma\alpha}\cdot a)(\gamma_{(\sigma\alpha)^{-1}}(c))
= (q_{\sigma\alpha}\circ \gamma_{(\sigma\alpha)^{-1}})(c)a
\equal{\equref{8.12.1}} q_e(c)a\\
&=& (q_e\cdot a)(c)
= (q_e\cdot a\circ \gamma_{(\sigma\alpha)^{-1}}^{-1})(\gamma_{(\sigma\alpha)^{-1}}(c)),
\end{eqnarray*}
as needed. This finishes our proof.
\end{proof}

\section{Galois group corings and graded Morita contexts}\selabel{9}
Let us call a $G$-$A$-coring $\ul{\Cc}=(\Cc_\alpha)_{\alpha\in G}$ a left homogeneous progenerator if every $\Cc_\alpha$ is a left $A$-progenerator.
We will now apply the Morita theory developed in the previous Section to find some equivalent properties for a left homogeneous progenerator group coring to be Galois.
Recall from \cite{C03} the following Theorem.

\begin{theorem}\thlabel{progen}
Let $(\Cc_e,x_e)$ be an $A$-coring with a fixed grouplike element, and assume that $\Cc_e$ is a left $A$-progenerator. We take a subring $B$ of $T_e=A^{{\rm co}\Cc_e}=\{a\in A~|~ax_e=x_ea\}$, and consider the map
$$\can_e':\ \Dd'=A\ot_{B}A\to \Cc_e,~ \can'(a\ot_B b)=ax b.$$
Then the following statements are equivalent:
\begin{enumerate}
\item \begin{itemize}
\item $\can_e'$ is an isomorphism of corings;
\item $A$ is faithfully flat as a left $B$-module.
\end{itemize}
\item \begin{itemize}
\item ${}^*\can_e'$ is an isomorphism of rings;
\item $A$ is a left $B$-progenerator.
\end{itemize}
\item \begin{itemize}
\item $B=T_e$;
\item the Morita context $\MM_e=(T_e,R_e,A,Q_e,\varphi_e,\psi_e)$ is strict.
\end{itemize}
\item \begin{itemize}
\item $B=T_e$;
\item $(F_8,G_8)$ is an equivalence of categories.
\end{itemize}
\end{enumerate}
\end{theorem}

Here $\MM_e$ is the Morita context introduced in the light of  \prref{8.9}, and $(F_8,G_8)$ the adjoint pair of functors considered before \leref{5.11}.
In the next Theorem we use the graded Morita context $\GG\MM$ of \thref{8.8}, and the adjoint pair of functors $(F_7,G_7)$ of \prref{5.4}.

\begin{theorem}
Let $(\ul{\Cc},\ul{x})$ be a left homogeneous progenerator $G$-$A$-coring with a fixed grouplike family. We take a subring $B$ of $T=A^{{\rm co}\ul{\Cc}}\subset T_e$, and consider the map
$$\ul{\can'}:\ \ul{\Dd'}=(A\ot_{B}A)\langle G\rangle \to \ul{\Cc},~ \can'_\alpha(\mu_\alpha(a\ot_B b))=ax_\alpha b.$$
Then the following statements are equivalent:
\begin{enumerate}
\item \begin{itemize}
\item $\ul{\can'}$ is an isomorphism of group corings;
\item $A$ is faithfully flat as a left $B$-module.
\end{itemize}
\item \begin{itemize}
\item ${}^*\ul{\can'}$ is an isomorphism of graded rings;
\item $A$ is a left $B$-progenerator.
\end{itemize}
\item \begin{itemize}
\item $B=T\cong S$;
\item the graded Morita context $\GG\MM=(G*S,R,A\{G\},QG,\omega,\nu)$ is strict.
\end{itemize}
\item \begin{itemize}
\item $B=T$;
\item $(F_7,G_7)$ is an equivalence of categories.
\end{itemize}
\end{enumerate}
\end{theorem}

\begin{proof}
$\ul{1)\Rightarrow 2)}$. Obviously ${}^*\ul{\can'}$ is an isomorphism if $\ul{\can'}$ is an isomorphism. In particular, $\can'_e$ is an isomorphism of corings, hence it follows from \thref{progen} that $A$ is a left $B$-progenerator. \\
$\ul{2)\Rightarrow 1)}$. Suppose that ${}^*\ul{\can'}:\ {}^*\ul{\Cc}=R\to {}^*\ul{\Dd'}=R'$ is an isomorphism. We then have that the right dual $({}^*\ul{\can'})^*:\ R'^*\to R^*,\ ({}^*\ul{\can'})^*(\varphi)=\varphi\circ {}^*\ul{\can'} $ is also an isomorphism. Since $\ul{\Cc}$ and $\ul{\Dd'}$ are left homogeneously finite, this map can be interpreted as the isomorphism
$$f=\iota^{-1}\circ ({}^*\ul{\can'})^*\circ \iota' :\prod_{\alpha\in G} \Dd'_\alpha \to \prod_{\alpha\in G}\Cc_\alpha,$$
where we denoted respectively $\iota$ and $\iota'$ for the isomorphisms $\prod_{\alpha\in G}\Cc_\alpha\cong R^*$ and $\prod_{\alpha\in G} \Dd'_\alpha \cong R'^*$ (see the beginning of \seref{4}). For all $\ul{d}=(d_\alpha)_{\alpha\in G}\in  \prod_{\alpha\in G} \Dd'_\alpha$ we have that
\begin{eqnarray*}
&&\hspace{-1.5cm}
f(\ul{d})= (\iota^{-1}\circ ({}^*\ul{\can'})^*\circ \iota')(\ul{d})= \iota^{-1}(\iota'(\ul{d})\circ {}^*\ul{\can'})\\
&=& \big( (\iota'(\ul{d})\circ {}^*\ul{\can'})(f^{(\alpha)})c^{(\alpha)} \big)_{\alpha\in G}= 
\big( \iota'(\ul{d})( f^{(\alpha)}\circ \can'_\alpha)c^{(\alpha)} \big)_{\alpha\in G}\\
&=& \big( f^{(\alpha)}(\can'_\alpha(d_\alpha))c^{(\alpha)} \big)_{\alpha\in G}= ( \can'_\alpha(d_\alpha))_{\alpha\in G},
\end{eqnarray*}
i.e., $f=\prod_{\alpha\in G}\can'_\alpha$. Now, since $f=\prod_{\alpha\in G}\can'_\alpha$ is an isomorphism, it follows that all $\can'_\alpha$ are isomorphisms. Indeed, $(\can'_\alpha)^{-1}=p'_\alpha\circ f^{-1}\circ i_\alpha:\Cc_\alpha \to \Dd'_\alpha$, where $i_\alpha$ and $p'_\alpha$ are the canonical injections and projections, respectively.
So the inverse of $\ul{\can'}$ is given by $(\ul{\can'})^{-1}=((\can'_\alpha)^{-1})_{\alpha\in G}$. 
Finally, since in particular ${}^*\can'_e$ is an isomorphism, \thref{progen} implies that $A$ is faithfully flat as a left $B$-module.\\
$\ul{1)\Rightarrow 3)}$. As in the proof of \prref{5.10}, $\ul{\Cc}$ is a cofree group coring with a fixed grouplike family, since $\ul{\can'}:\ \ul{\Dd'}\to \ul{\Cc}$ is an isomorphism,
and $\ul{\Dd'}$ is cofree. By \leref{5.9} we then have that $T=T_e$. Hence it follows from \thref{progen} and the fact that $\can'_e$ is an isomorphism that $B=T$ and that the Morita context $\MM_e$ is strict. It is easily verified that the graded Morita context $\MM_e[G]$ then also is strict, and likewise $\GG\MM$, see \thref{8.12}. From \prref{8.10} we get $T\cong S$.\\
$\ul{4)\Rightarrow 1)}$. By \prref{5.7} $\ul{\can}$ is an isomorphism, i.e. $(\ul{\Cc},\ul{x})$ is Galois. Hence (see \prref{5.10}) $(\ul{\Cc},\ul{x})$ is a cofree group coring with a fixed grouplike familie, and thus \leref{5.9} implies that $B=T=T_e$. Since $(F_7\cong F_2\circ F_8, G_7\cong G_8\circ G_2)$ and $(F_2,G_2)$ are equivalences (see \leref{5.11} and \thref{2.2}), it follows that also $(F_8,G_8)$ is an equivalence of categories. Finally it follows from \thref{progen} that $A$ is faithfully flat as a left $B$-module.\\
$\ul{3)\Rightarrow 4)}$. 
Suppose that $\GG\MM=(G*S,R,A\{G\},QG,\omega,\nu)$ is a strict graded Morita context. We then have a pair of inverse equivalences $(\widetilde{F}=-\ot_{G*S}A\{G\},\widetilde{G}=-\ot_R QG)$ between the categories $\Mm_{G*S}^G$ and $\Mm^G_R$. By \prref{8.10} we have that $G*S$ and $T[G]$ are isomorphic as graded rings. As a consequence the categories $\Mm^G_{G*S}$ and $\Mm^G_{T[G]}$ are isomorphic, and the latter is in turn isomorphic with $\Mm_T$ by the Structure Theorem for graded modules over strongly graded rings. Making use of the pair of functors $(F_3,G_3)$ which constitutes an isomorphism between $\Mm^{G,\ul{\Cc}}$ and $\Mm^G_R$ (see \prref{4.2}), we have the following pair $(\widetilde{F_7} ,\widetilde{G_7})$ of inverse equivalences between $\Mm_T$ and $\Mm^{G,\ul{\Cc}}$:
$$\xymatrix{
\widetilde{F_7}:\ \Mm_T\cong \Mm^G_{T[G]}\cong \Mm^G_{G*S}\ar@<.5ex>[r]^(.75){\widetilde{F}} & \Mm^G_R\ar@<.5ex>[l]^(.25){\widetilde{G}}\ar@<.5ex>[r]^(.4){G_3}&  \Mm^{G,\ul{\Cc}}\ar@<.5ex>[l]^(.6){F_3}\ :\widetilde{G_7}.
}$$
For $M\in \Mm_T$ we have that
$$\widetilde{F_7}(M)=\Big( \big( M[G] \ot_{G*S} A\{G\}\big)_\alpha \Big)_{\alpha\in G}\in \Mm^{G,\ul{\Cc}},
$$
where we denote $\big( M[G] \ot_{G*S} A\{G\}\big)_\alpha$ for the $\alpha$th homogeneous component of $M[G] \ot_{G*S} A\{G\} \in \Mm_R^G$.
The coaction maps of $\widetilde{F_7}(M)$ are given by
\begin{eqnarray*}
&&\hspace{-2cm}\widetilde{\rho}_{\alpha,\beta}:\ \big(M[G] \ot_{G*S} A\{G\}\big)_{\alpha\beta}\to \big(M[G] \ot_{G*S} A\{G\}\big)_\alpha\ot_A \Cc_\beta,\\
&&\hspace{-1cm}
\widetilde{\rho}_{\alpha,\beta}\Big(\sum_{\gamma\in G}m_\gamma u_\gamma \ot_{G*S} \mu_{\gamma^{-1}\alpha\beta}(a_\gamma)\Big)\\
&=&\sum_{\gamma\in G} m_\gamma u_\gamma\ot_{G*S} \mu_{\gamma^{-1}\alpha\beta}(a_\gamma)\cdot f^{(\beta)} \ot_A c^{(\beta)}\\
&=& \sum_{\gamma\in G} m_\gamma u_\gamma \ot_{G*S} \mu_{\gamma^{-1}\alpha}(f^{(\beta)}(x_\beta a_\gamma))\ot_A c^{(\beta)}.
\end{eqnarray*}
We now claim that $\widetilde{F_7}\cong F_7$. For $M\in \Mm_T$ and $\alpha\in G$, we consider the map
\begin{eqnarray*}
&&\hspace{-2cm}\varphi_{M,\alpha}:\ \big(M[G] \ot_{G*S} A\{G\}\big)_\alpha \to \mu_\alpha(M\ot_T A),\\
&&\hspace{-1cm}\varphi_{M,\alpha}\Big(\sum_{\gamma \in G} m_\gamma u_\gamma \ot_{G*S} \mu_{\gamma^{-1}\alpha}(a_\gamma)\Big)=\mu_\alpha\Big(\sum_{\gamma\in G}m_\gamma\ot_T a_\gamma\Big).
\end{eqnarray*}
 $\varphi_{M,\alpha}$ is well-defined, since
\begin{eqnarray*}
&&\hspace{-1.5cm}
\varphi_{M,\alpha}\Big(\sum_{\gamma \in G} m_\gamma u_\gamma\cdot u_e\ul{b} \ot_{G*S} \mu_{\gamma^{-1}\alpha}(a_\gamma)\Big)\\
&=&\varphi_{M,\alpha}\Big(\sum_{\gamma \in G} m_\gamma b_e u_\gamma \ot_{G*S} \mu_{\gamma^{-1}\alpha}(a_\gamma)\Big)
=\mu_\alpha\Big(\sum_{\gamma\in G}m_\gamma b_e \ot_T a_\gamma\Big)
\end{eqnarray*}
equals
\begin{eqnarray*}
&&\hspace{-1.5cm}
\varphi_{M,\alpha}\Big(\sum_{\gamma \in G} m_\gamma u_\gamma\ot_{G*S}  u_e\ul{b} \cdot \mu_{\gamma^{-1}\alpha}(a_\gamma)\Big)\\
&=& \varphi_{M,\alpha}\Big(\sum_{\gamma \in G} m_\gamma u_\gamma\ot_{G*S}   \mu_{\gamma^{-1}\alpha}(b_{\gamma^{-1}\alpha}a_\gamma)\Big)\\
&=&\varphi_{M,\alpha}\Big(\sum_{\gamma \in G} m_\gamma u_\gamma\ot_{G*S}  \mu_{\gamma^{-1}\alpha}(b_{e}a_\gamma)\Big)
=\mu_\alpha\Big(\sum_{\gamma\in G}m_\gamma  \ot_T b_ea_\gamma\Big),
\end{eqnarray*}
for all $\ul{b}=i(b_e)\in S$.
Clearly $\varphi_{M,\alpha}$ is right $A$-linear. Let us check that $\varphi_M=(\varphi_{M,\alpha})_{\alpha\in G}:\ \widetilde{F_7}(M)\to F_7(M)$ is a morphism in $\Mm^{G,\ul{\Cc}}$:
\begin{eqnarray*}
&&\hspace{-1.8cm}
((\varphi_{M,\alpha}\ot_A \Cc_\beta)\circ \widetilde{\rho}_{\alpha,\beta})\Big(\sum_{\gamma\in G} m_\gamma u_\gamma \ot_{G*S} \mu_{\gamma^{-1}\alpha\beta}(a_\gamma)\Big)\\
&=& 
\sum_{\gamma \in G}\varphi_{M,\alpha}\big(m_\gamma u_\gamma\ot_{G*S} \mu_{\gamma^{-1}\alpha}(f^{(\beta)}(x_\beta a_\gamma))\big) \ot_A c^{(\beta)}\\
&=&
\sum_{\gamma \in G}\mu_{\alpha}(m_\gamma \ot_T f^{(\beta)}(x_\beta a_\gamma)) \ot_A c^{(\beta)}\\
&=&
\sum_{\gamma \in G}\mu_{\alpha}(m_\gamma \ot_T 1_A) f^{(\beta)}(x_\beta a_\gamma) \ot_A c^{(\beta)}\\
&=&
\sum_{\gamma \in G}\mu_{\alpha}\big(m_\gamma \ot_T 1_A\big) \ot_A f^{(\beta)}(x_\beta a_\gamma) c^{(\beta)}\\
&=&
\sum_{\gamma \in G}\mu_{\alpha}\big(m_\gamma \ot_T 1_A\big) \ot_A x_\beta a_\gamma
=
\rho_{\alpha,\beta}\Big(\mu_{\alpha\beta}\Big(\sum_{\gamma\in G}m_\gamma\ot_T a_\gamma \Big)\Big)\\
&=& 
(\rho_{\alpha,\beta}\circ \varphi_{M,\alpha\beta})\Big( \sum_{\gamma\in G} m_\gamma u_\gamma \ot_{G*S} \mu_{\gamma^{-1}\alpha\beta}(a_\gamma)\Big).
\end{eqnarray*}
Let us finally show that $\varphi_{M}$ is an isomorphism in $\Mm^{G,\ul{\Cc}}$. It suf\emph{}fices to check that the inverse of $\varphi_{M,\alpha}$ is given by 
$$\varphi_{M,\alpha}^{-1}\Big( \mu_\alpha\Big( \sum_{i=1}^n m_i\ot_T a_i \Big)\Big)= \sum_{i=1}^n m_iu_e \ot_{G*S} \mu_\alpha(a_i)= \sum_{i=1}^n m_iu_\alpha \ot_{G*S} \mu_e(a_i):$$
\begin{eqnarray*}
&&\hspace{-1.3cm}
(\varphi_{M,\alpha}^{-1}\circ \varphi_{M,\alpha})\Big(\sum_{\gamma\in G} m_\gamma u_\gamma \ot_{G*S} \mu_{\gamma^{-1}\alpha}(a_\gamma)\Big)\\
&=&
\varphi_{M,\alpha}^{-1}\Big(\mu_\alpha\Big(\sum_{\gamma\in G}m_\gamma\ot_T a_\gamma\Big)\Big)
= \sum_{\gamma\in G} m_\gamma u_e \ot_{G*S} \mu_\alpha(a_\gamma)\\
&=&
\sum_{\gamma\in G} m_\gamma u_e \ot_{G*S} \mu_{\gamma \gamma^{-1}\alpha}(a_\gamma)
=\sum_{\gamma\in G} m_\gamma u_e \ot_{G*S} u_\gamma\ul{1}\cdot  \mu_{ \gamma^{-1}\alpha}(a_\gamma)\\
&=&
\sum_{\gamma\in G} m_\gamma u_e \cdot u_\gamma\ul{1}\ot_{G*S}  \mu_{ \gamma^{-1}\alpha}(a_\gamma)
= \sum_{\gamma\in G} m_\gamma u_\gamma \ot_{G*S} \mu_{\gamma^{-1}\alpha}(a_\gamma);
\end{eqnarray*}
\begin{eqnarray*}
&&\hspace{-1.3cm}
(\varphi_{M,\alpha}\circ \varphi_{M,\alpha}^{-1})\Big( \mu_\alpha\Big( \sum_{i=1}^n m_i\ot_T a_i \Big)\Big)
= \varphi_{M,\alpha}\Big(\sum_{i=1}^n m_iu_e \ot_{G*S} \mu_\alpha(a_i)\Big)\\
&=& \sum_{i=1}^n \mu_\alpha(m_i\ot_T a_i)
= \mu_\alpha\Big( \sum_{i=1}^n m_i\ot_T a_i \Big).
\end{eqnarray*}
So we have shown that $F_7$ and $\widetilde{F_7}$ are naturally isomorphic. From the uniqueness of the adjoint functor, it follows that also $G_7 \cong \widetilde{G_7}$. The fact that $(\widetilde{F_7},\widetilde{G_7})$ is a pair of inverse equivalences implies that also $(F_7,G_7)$ is a pair of inverse equivalences, as needed.
\end{proof}

\section{Application to $\ul{H}$-comodule algebras}\selabel{10}
Let $k$ be a  commutative ring.
Recall \cite{Turaev} that a Hopf $G$-coalgebra is a $G$-coalgebra $\ul{H}=(H_\alpha)_{\alpha\in G}$
with the following additional structure: every $H_\alpha$ is a $k$-algebra, such that
$\Delta_{\alpha,\beta}$ and $\epsilon$ are algebra maps, and we have a family of maps
$S_\alpha:\ H_{\alpha^{-1}}\to H_\alpha$ such that
$$S_\alpha(h_{(1,\alpha^{-1})})h_{(2,\alpha)}=h_{(1,\alpha)}S_\alpha(h_{(2,\alpha^{-1})})=\varepsilon(h)1_{H_\alpha},$$
for every $h\in H_e$. A right $G$-$\ul{H}$-comodule algebra is a $k$-algebra $A$ with
a right $\ul{H}$-coaction $\ul{\rho}=(\rho_\alpha)_{\alpha\in G}$ such that
$$\rho_\alpha(ab)=a_{[0]}b_{[0]}\ot a_{[1,\alpha]}b_{[1,\alpha]}~~{\rm and}~~
\rho_\alpha(1_A)=1_A\ot 1_{H_\alpha},$$
for all $a,b\in A$ and $\alpha\in G$. This notion was introduced by the third author in \cite{Wang1}.
The proof of the following result is straightforward.

\begin{proposition}\prlabel{10.1}
Let $\ul{H}$ be a Hopf $G$-coalgebra, and $A$ right $G$-$\ul{H}$-comodule algebra.
Then $\ul{\Cc}=A\ot \ul{H}=(A\ot H_\alpha)_{\alpha\in G}$ is a $G$-$A$-coring. The $A$-bimodule
structures are given by the formulas
$$a'(b\ot h)a=a'ba_{[0]}\ot ha_{[1,\alpha]},$$
for all $a,a',b\in A$, $\alpha\in G$, $h\in H_\alpha$. The comultiplication and counit maps are the following:
$$\Delta_{\alpha,\beta}:\ A\ot H_{\alpha\beta}\to (A\ot H_\alpha)\ot_A(A\ot H_\beta),~~
\Delta_{\alpha,\beta}(a\ot h)=(a\ot h_{(1,\alpha)})\ot_A (1_A\ot h_{(2,\beta)});$$
$$\varepsilon=A\ot\epsilon:\ A\ot H_e\to A.$$
$\ul{x}=(1_A\ot 1_{H_\alpha})_{\alpha\in G}$ is a grouplike family of $A\ot\ul{H}$.
\end{proposition}

It is easy to see that, for $a\in A$, $a\in A^{{\rm co}\ul{\Cc}}$ if and only if
$a\ot 1_{H_\alpha}$ equals $(1_A\ot 1_{H_\alpha})a=a_{[0]}\ot a_{[1,\alpha]}=
\rho_\alpha(a)$, for all $\alpha\in G$. With notation as in \cite{Wang4}, this means that
$$A^{{\rm co}\ul{\Cc}}=A^0=\{a\in A~|~\rho_\alpha(a)=a\ot 1_{H_\alpha},~{\rm for~all}~
\alpha\in G\}.$$
Let $B\to A^{{\rm co}\ul{\Cc}}$ be a ring morphism.
We can then compute the morphism $\ul{\can}:\ (A\ot_{B}A)\lan G\ran \to
A\ot \ul{H}$ as follows: $\can_\alpha:\ \mu_\alpha(A\ot_B A)\to A\ot H_\alpha$ is given by the formula
$$\can_\alpha(\mu_\alpha(a\ot b))=a(1_A\ot 1_{H_\alpha})b=ab_{[0]}\ot b_{[1,\alpha]}.$$
This proves the following result.

\begin{proposition}\prlabel{10.2}
Let $A$ be a right comodule algebra over a Hopf $G$-coalgebra $\ul{H}$.
Then $(A\ot \ul{H},(1_A\ot 1_{H_\alpha})_{\alpha\in G})$ is a Galois $G$-$A$-coring
if and only if $A$ is a $G$-$\ul{H}$-Galois extension of $A^{{\rm co}\ul{\Cc}}=A^0$,
in the sense of \cite[Def. 7.1]{Wang4}.
\end{proposition}

Let $H$ be a Hopf algebra, and $\ul{H}=(H_\alpha)_{\alpha\in G}$ a set of isomorphic
copies of $H$, indexed by the group $G$. Let $H_e=H$, and $\lambda_\alpha:\
H\to H_\alpha$ the connecting isomorphism. Then $\ul{H}$ is a Hopf $G$-coalgebra,
with structure maps
$$\Delta_{\alpha,\beta}(\lambda_{\alpha\beta}(h))=\lambda_{\alpha}(h_{(1)})
\ot \lambda_{\beta}(h_{(2)});$$
$$S_\alpha(\lambda_{\alpha^{-1}}(h))=\lambda_\alpha(S(h)).$$
The counit is the counit of $H$, and every $H_\alpha$ is a $k$-algebra. We call
$\ul{H}=H\lan G\ran $ the cofree Hopf $G$-coalgebra associated to $H$. Using 
Propositions \ref{pr:5.10} and \ref{pr:10.2}, we obtain the following result:

\begin{proposition}\prlabel{10.3}
Let $A$ be a right comodule algebra over a Hopf $G$-coalgebra $\ul{H}$.
$A$ is a $G$-$\ul{H}$-Galois extension of $A^0$ if and only if
$\ul{H}$ is a cofree Hopf $G$-coalgebra and $A$ is an $H$-Galois extension
of $A^{{\rm co}H_e}=A^0$.
\end{proposition}

A right relative $(\ul{H},A)$-Hopf module (in \cite{Wang4} termed a right
$G$-$(\ul{H},A)$-Hopf module) is a right $A$-module $M$, with the additional
structure of right $\ul{H}$-comodule, such that the compatibility condition
$$\rho_\alpha(ma)=m_{[0]}a_{[0]}\ot m_{[1,\alpha]}a_{[1,\alpha]}$$
holds for all $m\in M$, $a\in A$, $\alpha\in G$. $\Mm_A^{\ul H}$ will denote the category
of right relative $(\ul{H},A)$-Hopf modules.\\
In a similar way, a right relative group $(\ul{H},A)$-Hopf module is a family of right
$A$-modules $(M_\alpha)_{\alpha\in G}$, with the additional structure 
$(\rho_{\alpha,\beta})_{\alpha,\beta\in G}$ of right $G$-$\ul{H}$-comodule, with
the compatibility relation
$$\rho_{\alpha\beta}(ma)=m_{[0,\alpha]}a_{[0]}\ot m_{[1,\alpha]}a_{[1,\beta]},$$
for all $\alpha,\beta\in G$, $m\in M_{\alpha\beta}$ and $a\in A$. The category of
right relative group $(\ul{H},A)$-Hopf modules is denoted by $\Mm_A^{G,{\ul H}}$.
The proof of the following result is straightforward, and is left to the reader.

\begin{proposition}\prlabel{10.4}
Let $A$ be a right comodule algebra over a Hopf $G$-coalgebra $\ul{H}$. Then
we have isomorphisms of categories
$\Mm_A^{\ul H}\cong \Mm^{A\ot\ul{H}}$ and $\Mm_A^{G,{\ul H}}\cong \Mm^{G,A\ot\ul{H}}$.
\end{proposition}

Let $B\to A^0$ be a ring morphism.
It follows from Propositions \ref{pr:5.4} and \ref{pr:10.4} that we have a pair of adjoint
functors $(F_7,G_7)$ between $\Mm_B$ and $\Mm_A^{G,{\ul H}}$. As an application
of \thref{5.12}, we obtain the following Structure Theorem for relative group $(\ul{H},A)$-Hopf modules.

\begin{proposition}\prlabel{10.5}
Let $A$ be a right comodule algebra over a Hopf $G$-coalgebra $\ul{H}$,
and $B\to A^0$ a ring morphism. Then
the following assertions are equivalent.
\begin{enumerate}
\item $B\cong A^0$, $A$ is a $G$-$\ul{H}$-Galois extension of $A^0$, and
$A$ is faithfully flat as a left $B$-module;
\item $(F_7,G_7)$ is a pair of inverse equivalences and 
$A$ is flat as a left $B$-module.
\end{enumerate}
\end{proposition}

Let us finally compute the left dual graded $A$-ring $R$ of $A\ot \ul{H}$. We have an
isomorphism of $k$-modules
$$R=\bigoplus_{\alpha\in G} {}_A\Hom(A\ot H_{\alpha^{-1}},A)\cong
\bigoplus_{\alpha\in G} \Hom(H_{\alpha^{-1}},A).$$
The multiplication (and the $A$-bimodule structure) on $R$ can be transported to
$\bigoplus_{\alpha\in G} \Hom(H_{\alpha^{-1}},A)$. We obtain the following multiplication
rule, for $f\in \Hom(H_{\alpha^{-1}},A)\cong R_\alpha$, 
$g\in \Hom(H_{\beta^{-1}},A)\cong R_\beta$, $h\in H_{(\alpha\beta)^{-1}}$:
\begin{equation}\eqlabel{10.6.1}
(f\#g)(h)=f(h_{(2,\alpha^{-1})})_{[0]} g\bigl(h_{(1,\beta^{-1})}f(h_{(2,\alpha^{-1})})_{[1,\beta^{-1}]}\bigr).
\end{equation}

Before we investigate more carefully the situation where $\ul{H}$ is homogeneously finite
(that is, every $H_\alpha$ is a finitely generated and projective $k$-module, we make
some general observations.\\

Let $K$ be a (classical) Hopf algebra, and $A$ a left $K$-module algebra. Then we
can form the smash product $K^{\rm op}\# A$, with multiplication rule
\begin{equation}\eqlabel{10.6.2}
(h\#a)(k\#b)=k_{(1)}h\# (k_{(2)}\cdot a)b.
\end{equation}
It is well-known that $K^{\rm op}\# A$ is an $A$-ring.\\
We call $K$ a graded Hopf algebra if $K$ is a Hopf algebra and a $G$-graded algebra
such that $\Delta(K_\alpha)\subset K_\alpha\ot K_\alpha$ and $S(K_\alpha)
\subset K_{\alpha^{-1}}$. This implies in particular that every $K_\alpha$ is a subcoalgebra
of $K$. If $K$ is a graded Hopf algebra, and $A$ is a left $K$-module algebra, then
$K^{\rm op}\# A$ is a graded $A$-ring.\\
In \cite{Zunino1}, a $G$-graded Hopf algebra is called a Hopf $G$-algebra in packed form.
The defining axioms of a Hopf $G$-algebra are formally dual to the defining axioms of
a Hopf $G$-coalgebra. A Hopf $G$-algebra is a family of $k$-coalgebras $\ul{K}=
(K_\alpha)_{\alpha\in G}$ together with $k$-coalgebra maps
$\mu_{\alpha,\beta}:\ K_\alpha\ot K_\beta\to K_{\alpha\beta}$ and $\eta:\ k\to K_e$
satisfying the obvious associativity and unit properties. We also need antipode maps
$S_\alpha:\ K_\alpha\to K_{\alpha^{-1}}$ such that 
$$\mu_{\alpha^{-1},\alpha}(S_\alpha(k_{(1)})\ot k_{(2)})=
\mu_{\alpha,\alpha^{-1}}(k_{(1)}\ot S_\alpha(k_{(2)}))=\eta(\varepsilon(k)),$$
for all $k\in K_\alpha$. It is straightforward to show that $K=\bigoplus_{\alpha\in G}K_\alpha$
is a graded Hopf algebra. Conversely, if $K$ is a graded Hopf algebra, then
$(K_\alpha)_{\alpha\in G}$ is a Hopf $G$-algebra. Thus we have an isomorphism between
the categories of $G$-graded Hopf algebras and Hopf $G$-algebras.\\
If $\ul{H}$ is a homogeneously finite Hopf $G$-coalgebra, then $\ul{K}=
(H^*_{\alpha^{-1}})_{\alpha\in G}$ is a Hopf $G$-algebra, and, consequently,
$K=\bigoplus_{\alpha\in G} H^*_{\alpha^{-1}}$ is a $G$-graded Hopf algebra.\\
If $A$ is a right $\ul{H}$-module algebra, then it is also a left $K$-module algebra, with
action $h^*\cdot a=\lan h^*,a_{[1,\alpha^{-1}]}\ran a_{[0]}$, for all $h^*\in K_\alpha=
H^*_{\alpha^{-1}}$.\\
For every $\alpha\in G$,
$${}_A\Hom(A\ot H_{\alpha^{-1}},A)\cong \Hom(H_{\alpha^{-1}},A)\cong
H^*_{\alpha^{-1}}\ot A$$
is the degree $\alpha$ component of $K^{\rm op}\# A$.

\begin{theorem}\thlabel{10.6}
Let $\ul{H}$ be a homogeneously finite Hopf $G$-coalgebra, and $A$ a right
$\ul{H}$-comodule algebra. Then $R=\bigoplus_{\alpha\in G} {}_A\Hom(A\ot H_{\alpha^{-1}},A)$
is isomorphic to $K^{\rm op}\# A$ as a $G$-graded $A$-ring. Consequently, the
categories $\Mm_A^{G,\ul{\Cc}}$ and $\Mm_{K^{\rm op}\# A}^G$ are isomorphic.
\end{theorem}

\begin{proof}
We have to show that the $k$-module isomorphisms
$$\lambda_\alpha:\ H^*_{\alpha^{-1}}\ot A\to \Hom(H_{\alpha^{-1}},A),~~
\lambda_\alpha(h^*\ot a)(h)=\lan h^*,h\ran a$$
transport the multiplication rule \equref{10.6.2} to \equref{10.6.1}. Take
$\alpha,\beta\in G$, $h^*\in H^*_{\alpha^{-1}}$, $k^*\in H^*_{\beta^{-1}}$, $a,b\in A$,
and write $f=\lambda_\alpha(h^*\ot a)$, $g= \lambda_\beta(k^*\ot b)$. For
$h\in H_{(\alpha\beta)^{-1}}$, we have
\begin{eqnarray*}
&&\hspace*{-2cm}
(f\#g)(h)=
\lan h^*,h_{(2,\alpha^{-1})}\ran a_{[0]} g(h_{(1,\beta^{-1})}a_{[1,\beta^{-1}]})\\
&=&
\lan h^*,h_{(2,\alpha^{-1})}\ran \lan k^*, h_{(1,\beta^{-1})}a_{[1,\beta^{-1}]}\ran a_{[0]}b\\
&=&\lan h^*,h_{(2,\alpha^{-1})}\ran   \lan k_{(1)}^*, h_{(1,\beta^{-1})}\ran 
\lan k_{(2)}^*,a_{[1,\beta^{-1}]}\ran a_{[0]}b\\
&=& \lan k_{(1)}^**h^*,h\ran (k_{(2)}^*\cdot a)b\\
&=& \lambda_{\alpha\beta}(k_{(1)}^**h^*\# (k_{(2)}^*\cdot a)b)(h),
\end{eqnarray*}
and we conclude that
$$\lambda_\alpha(h^*\ot a)\# \lambda_\beta(k^*\ot b)=\lambda_{\alpha\beta}
((h^*\# a)(k^*\# b)),$$
as needed.
\end{proof}

\end{document}